\documentclass[11pt]{article}
\usepackage{amsmath}
\usepackage{amssymb}
\usepackage{amscd}
\usepackage[all]{xy}
\usepackage[T1]{fontenc}
\usepackage{epsfig}
\newtheorem{theo}{Theorem}[section]
\newtheorem{prop}[theo]{Proposition}
\newtheorem{coro}[theo]{Corollary}
\newtheorem{lemm}[theo]{Lemma}
\newtheorem{defi}[theo]{Definition}

\newtheorem{ex}[theo]{Example}

\title{\textbf{\textsc{galois and bigalois objects over monomial non 
semisimple hopf algebras}}} 
\author{\textsc{Julien Bichon}}
\date{\small{\textsl{Laboratoire de Math\'ematiques Appliqu\'ees,
Universit\'e de Pau et des Pays de l'Adour, \\
IPRA, Avenue de l'universit\'e,
64000 Pau, France.}
E-mail: Julien.Bichon@univ-pau.fr}}

\makeatletter
\renewcommand{\@makefnmark}{}
\makeatother

\begin{document}

\maketitle

\begin{abstract}
We describe the Hopf-Galois extensions of the base field
and the biGalois groups
of non semisimple monomial Hopf algebras. 
The main feature of our description is the use
of modified versions of the second cohomology group
of the grouplike elements.
Our computations generalize the 
previous ones of Masuoka and Schauenburg for the Taft algebras.
\end{abstract}


\section*{Introduction}

This paper contributes to the general problem
of determining the Hopf-Galois extensions of a given Hopf algebra.
We describe here the Hopf-Galois extensions (of the base field)
and the biGalois groups of an important class of pointed Hopf algebras, namely
the monomial non semisimple Hopf algebras, which were recently classified
by Chen, Huang, Ye and Zhang \cite{[CHYZ]}.

The first classification results of Hopf-Galois extensions
for a non commutative and non cocommutative Hopf algebra
are due to Masuoka \cite{[Ma]}, 
without any restriction on the subring of invariants,
in the case of the Taft algebras. 
Then Schauenburg \cite{[Sc2]} determined all Hopf-Galois extensions
of the base field (= Galois objects) for the general Taft algebras with several
grouplike elements, and computed their biGalois groups. The 
generalized Taft algebras are monomial Hopf algebras, 
so our results generalize Schauenburg's ones. 
The main feature of our description of Galois objects and biGalois
groups for monomial Hopf algebras is the use of modified versions
of the second cohomology group of the grouplike elements.
It gives a general concise descriptive formula for the set
of isomorphism classes of Galois objects and for the biGalois
groups, leaving the explicit computations to be done
at the cohomology group level. 

Our computations do not use the underlying quiver of a monomial Hopf
algebra: it certainly would be interesting to get a connection
between our description and the properties of the quiver.
Also it would interesting to describe Galois objects and biGalois groups
over the Hopf algebra associated to a Hopf quiver 
by Cibils and Rosso in \cite{[CR]}.

\medskip

This paper is organized as follows.
The main definitions and results around Hopf-Galois extensions
and  biGalois groups are recalled in the first section.
We also recall the notion of a group datum of \cite{[CHYZ]},
which by the results of \cite{[CHYZ]} 
is equivalent to the notion of a monomial Hopf algebra, 
and divide the various group data into six different types. 
Finally this preliminary section is concluded with some
cohomogical preliminaries: the definition of the
modified second cohomology group, to be used
in the description results of the next sections.
In Section 2, we describe the Galois objects over the 
(monomial) Hopf algebra associated with a group datum, the type VI
case being treated thanks to the results of the third section.
The very end of the section deals with the classification
up to homotopy of these Galois extensions, a concept
recently introduced by Kassel and Schneider \cite{[KS]}: 
this is an easy task using the results of \cite{[KS]}.
Section 3 contains the general description of the BiGalois group
of a monomial Hopf algebra.
In Section 4, some explicit examples are examined in detail:
Hopf algebras associated with cyclic group data, including
the Taft algebras and the simple-pointed Hopf algebras,
and Hopf algebras associated with decomposable group data, including
the generalized Taft algebras with several grouplike.
In particular the results of Masuoka \cite{[Ma]} and Schauenburg \cite{[Sc2]}
are recovered.

\medskip

Throughout this paper $k$ is a commutative field containing all
primitive roots of unity. This forces $k$ to be a characteristic
zero field. The multiplicative group of non-zero elements
of $k$ is denoted by $\dot{k}$. For $n \in \mathbb N^*$, the
group of nth-roots of unity in $k$ is denoted
by $\mu_n$. By our assumption we have $|\mu_n| = n$.

\section{Preliminaries}

\subsection{Hopf-Galois extensions}

We first recall the principal facts concerning
Hopf-Galois extensions,
the main object of study in this paper.
The book \cite{[Mo]} is a convenient reference for this topic. 

Let $A$ be a Hopf algebra.
A \textbf{right $A$-Galois extension (of $k$)}
is a non-zero right $A$-comodule algebra $Z$ such that the linear map
$\kappa_r$ defined by the composition
\begin{equation*}
\begin{CD}
\kappa_r : Z \otimes Z @>1_Z \otimes \alpha>>
Z \otimes Z \otimes A @>m_Z \otimes 1_A>> Z \otimes A
\end{CD}
\end{equation*}
where $\alpha$ is the coaction of $A$
and $m_Z$ is the multiplication of $Z$, is bijective.
We also say that a right $A$-Galois extension (of $k$)
is an \textbf{$A$-Galois object}.
A morphism of $A$-Galois objects is an $A$-colinear algebra 
morphism. 
It is known that any morphism of $A$-Galois extensions is an isomorphism.
The set of isomorphism classes of $A$-Galois objects
is denoted by Gal$(A)$. For example, 
when $A = k[G]$ is a group algebra, 
we have ${\rm Gal}(k[G]) \cong H^2(G,\dot k)$.
However Gal$(A)$ does not carry a natural group structure in general.

Amongst $A$-Galois extensions, the cleft ones will be especially important
in this paper. First recall 
that a 2-cocycle over $A$ is a convolution invertible linear map
$\sigma : A \otimes A \longrightarrow k$ satisfying
$$\sigma(a_{(1)}, b_{(1)}) \sigma(a_{(2)}b_{(2)},c) =
\sigma(b_{(1)},c_{(1)}) \sigma(a,b_{(2)} c_{(2)})$$
and $\sigma(a,1) = \sigma(1,a) = \varepsilon(a)$, for all $a,b,c \in A$.
For example, when $A = k[G]$ is a group algebra, 
a cocycle over $A$ corresponds precisely to an element of $Z^2(G, \dot k)$.
To any 2-cocycle $\sigma$ on $A$ we associate the right $A$-comodule
algebra $_\sigma \! A$. As an $A$-comodule
 $_{\sigma} \! A = A$ and the product
of $_{\sigma}A$ is defined to be
$$a _{\sigma} \! . b = \sigma(a_{(1)}, b_{(1)}) a_{(2)} b_{(2)}, 
\quad a,b \in A.$$
Let $Z$ be a right $A$-comodule algebra. Then $Z$ is said to be \textbf{cleft}
if one the following equivalent conditions holds:

\noindent
$\bullet$ There exists a 2-cocycle over $A$ such that
$Z \cong {_\sigma \!A}$ as $A$-comodule algebras.

\noindent
$\bullet$ $Z$ is $A$-Galois and $Z \cong A$ as $A$-comodules.

\noindent
$\bullet$ There exists an $A$-colinear convolution invertible
map $A \longrightarrow Z$. 

\medskip

The equivalence of the above assertions are due to Doi and Takeuchi
and to Blattner and Montgomery.
We refer the reader to \cite{[Mo]} for proofs and for the original references.
It is known that any $A$-Galois extension of $k$ is cleft if $A$ has
one of the following properties:

\noindent
$\bullet$ $A$ is finite-dimensional (by \cite{[KC]}),
 
\noindent 
$\bullet$ $A$ is pointed, i.e. all its simple comodules are
one-dimensional, (by  \cite{[Gu]}, Remark 10).

Since we only consider pointed (and finite-dimensional) 
Hopf algebras in this paper, we only consider cleft Galois objects.
The reader might consult \cite{[Bi]} for examples of non cleft
Galois objects.

\medskip

Although this will not be needed in this paper, the reader
might like to know about the categorical interpretation
of Hopf-Galois extensions, which is especially enlightening
when we have a view towards biGalois objects.
Let Comod$(A)$ be the category of left $A$-comodules.
Recall \cite{[U]} that a fibre functor 
$\omega : {\rm Comod}(A) \longrightarrow {\rm Vect}(k)$ 
is a monoidal $k$-linear exact and faithful functor commuting with colimits.
Let $Z$ be a right $A$-Galois extension of $k$: Ulbrich \cite{[U]}
associates a fibre functor $\omega_Z$ to such a Galois extension.
This fibre functor is defined by 
$\omega_Z(V) = Z \square_A V$ where $Z\square_A V$ is the cotensor product
over $A$ of the right $A$-comodule $Z$ and of the left $A$-comodule $V$.
Then Ulbrich shows in \cite{[U]} that this defines a category equivalence
between $A$-Galois extensions of $k$ and fibre functors over Comod$(A)$.
The interpretation of the notion of cleft extension
is particularly nice in this setting.
Let $\omega$ be a fibre functor on Comod$(A)$. Then
the following assertions are equivalent:

\noindent
$\bullet$ The Galois extension corresponding to the
fibre functor $\omega$ is cleft.
 
\noindent
$\bullet$ The fibre functor $\omega$ is isomorphic, as a functor,
with the forgetful functor. 

\noindent
$\bullet$ The fibre functor $\omega$, when restricted to finite-dimensional
$A$-comodules, preserves the dimensions of the underlying vector spaces.

The hard part when showing the equivalence of the above
assertions is to see that the third one implies the second one.
This a consequence of a result of Etingof and Gelaki
(Proposition 4.2.2 in \cite{[EG]}) combined with Tannaka-Krein type
results.  

\medskip 

Let us now come to biGalois extensions.
Similarly to the right case,
a \textbf{left $A$-Galois extension} (of $k$)
is a non-zero left $A$-comodule algebra $Z$ such that the linear map
$\kappa_l$ defined by the composition
\begin{equation*}
\begin{CD}
\kappa_l : Z \otimes Z @>\beta \otimes 1_Z>>
A \otimes Z \otimes Z @>1_A \otimes m_Z>> A \otimes Z
\end{CD}
\end{equation*}
where $\beta$ is the coaction of $A$ and $m_Z$ is the multiplication
of $Z$, is bijective.

Let $A$ and $B$ be Hopf algebras. An algebra $Z$ 
is said to be \textbf{an $A$-$B$-bigalois extension} \cite{[Sc1]} if $Z$
is both a left $A$-Galois extension and a right $B$-Galois extension,
and if $Z$ is an $A$-$B$-bicomodule.

Now consider a Hopf algebra $A$ and an $A$-Galois object 
$Z$. Then Schauenburg shows (Theorem 3.5 in \cite{[Sc1]})
that there exists a Hopf algebra $L(Z,A)$ such that
$Z$ is an $L(Z,A)$-$A$-biGalois extension. The Hopf
algebra $L(Z,A)$ is unique: if $B$ is another Hopf algebra
such that $Z$ is $B$-$A$-biGalois, there exists a unique
Hopf algebra isomorphism $f : L(Z,A) \rightarrow B$ such that
$(f \otimes {\rm id}_Z) \circ \beta = \beta'$, where
$\beta$ and $\beta '$ denote the respective left coactions
of $L(Z,A)$ and $B$. In fact the Hopf algebra
$L(Z,A)$ is the Tannaka-Krein reconstructed object from
the fibre functor associated with the $A$-Galois extension
$Z$: this is Theorem 5.5 in \cite{[Sc1]}. 
In particular the $k$-linear monoidal comodule categories
over two Hopf algebras $A$ and $B$ are equivalent if and only if
there exists an $A$-$B$-biGalois extension. More
precisely by Corollary 5.7 of \cite{[Sc1]}
there is a category equivalence between:

\noindent
$\bullet$ The category of $k$-linear monoidal equivalences
of categories between Comod$(A)$ and Comod$(B)$. 

\noindent
$\bullet$ The category of $A$-$B$ biGalois extensions.

Thus BiGal$(A)$, the \textbf{biGalois group} of $A$, 
defined by Schauenburg in \cite{[Sc1]}
as the set of isomorphism classes of $A$-biGalois extensions
(i.e. $A$-$A$-biGalois extensions) endowed with the cotensor product,
might also be seen as the group of monoidal isomorphism
classes of  monoidal $k$-linear auto-equivalences of the category
Comod$(A)$. In this way if there exists an $A$-$B$-biGalois
extension, then there is a bijection 
Gal$(A) \cong {\rm Gal}(B)$ and a group isomorphism
BiGal$(A) \cong {\rm BiGal}(B)$.   
For example, when $A = k[G]$ is a group algebra,
we have a group isomorphism BiGal$(k[G]) \cong
{\rm Aut}(G) \ltimes H^2(G,\dot k)$.

\subsection{Group data and monomial Hopf algebras}

Following Chen, Huang, Ye and Zhang \cite{[CHYZ]}, we 
define a \textbf{group datum} (over $k$) to be
a quadruplet $\mathbb G = (G,g,\chi,\mu)$ consisting of
a finite group $G$, a central element $g \in G$,
a character $\chi : G \longrightarrow \dot k$ 
with $\chi(g) \not = 1$ and an element $\mu \in k$ such that
$\mu =0$ if $o(g) =o(\chi(g))$, and that if $\mu \not =0$,
then $\chi^{o(\chi(g))}=1$.

\smallskip

We need to classify the various group data
into six different types.  

\noindent
$\bullet$ A \textbf{type I group datum} $\mathbb G = (G,g,\chi,\mu)$
is a group datum with $\mu =0$, $d= o(\chi(g)) = o(g)$
and $\chi^d =1$. 
In this case we simply write $\mathbb G = (G,g,\chi)$.

\noindent
$\bullet$ A \textbf{type II group datum} $\mathbb G = (G,g,\chi,\mu)$
is a group datum with $\mu =0$, $d= o(\chi(g)) = o(g)$
and $\chi^d \not = 1$. 
In this case we simply write $\mathbb G = (G,g,\chi)$.

\noindent
$\bullet$ A \textbf{type III group datum} $\mathbb G = (G,g,\chi,\mu)$
is a group datum with $\mu =0$, $d= o(\chi(g)) < o(g)$
and $\chi^d =1$. 
In this case we simply write $\mathbb G = (G,g,\chi)$.

\noindent
$\bullet$ A \textbf{type IV group datum} $\mathbb G = (G,g,\chi,\mu)$
is a group datum with $\mu =0$, with $d= o(\chi(g)) < o(g)$,
with $\chi^d \not = 1$, and such that there does not exist
$\sigma \in Z^2(G, \dot k)$ with
$\chi^d(h) \sigma(h,g^d) = \sigma(g^d,h)$, $\forall h \in G$. 
In this case we simply write $\mathbb G = (G,g,\chi)$.

\noindent
$\bullet$ A \textbf{type V group datum} $\mathbb G = (G,g,\chi,\mu)$
is a group datum with $\mu =0$, $d= o(\chi(g)) < o(g)$
with $\chi^d \not = 1$, and such that there exists
$\sigma \in Z^2(G, \dot k)$ with
$\chi^d(h) \sigma(h,g^d) = \sigma(g^d,h)$, $\forall h \in G$. 
In this case we simply write $\mathbb G = (G,g,\chi)$.

\noindent
$\bullet$ A \textbf{type VI group datum} $\mathbb G = (G,g,\chi,\mu)$
is a group datum with $\mu \not =0$ (and hence 
$d = o(\chi(g)) < o(g)$ and $\chi^{o(\chi(g))}=1$). 

\medskip

Let $\mathbb G = (G,g,\chi,\mu)$ be a group datum. 
A Hopf algebra $A(\mathbb G)$ is associated with $\mathbb G$ 
in \cite{[CHYZ]}. We will slightly change the conventions
of \cite{[CHYZ]} for the formula defining the coproduct,
but this will not change the whole set of isomorphism classes.
As an algebra $A(\mathbb G)$ is the
quotient of the free product algebra
$k[x] * k[G]$ by the 
two-sided ideal generated by the relations
$$xh = \chi(h) hx, \ \forall h \in H , \quad \quad 
x^d = \mu(1-g^d), \ {\rm where} \ d = o(\chi(g)).$$  
The Hopf algebra structure of $A(\mathbb G)$ is defined
by:
$$\Delta(x) = 1  \otimes x + x \otimes g, \quad \varepsilon(x) =0, \quad
S(x) = -xg^{-1},$$
$$\Delta(h)= h \otimes h \ , \quad \varepsilon(h) =1 \ , \quad
S(h)=h^{-1} \ , \ \forall h \in G.$$ 
Using the diamond lemma \cite{[Be]}, it is not difficult
to see that the set 
$\{hx^i, 0 \leq i \leq d-1, h \in G \}$ is a linear basis
of $A(\mathbb G)$, and hence $\dim_k(A(\mathbb G))= |G|d$.

The Hopf algebras $A(\mathbb G)$ have been
shown in \cite{[CHYZ]} to be exactly the monomial non semisimple
Hopf algebras: see \cite{[CHYZ]} for the precise concept
of a monomial Hopf algebra. 

\medskip

Several particular cases of this Hopf algebra  
construction were considered 
in the literature: the Taft algebras \cite{[Ta]}, the simple-pointed
finite dimensional Hopf algebras (see \cite{[Ra]})...
See  section 4.
In fact a group datum $\mathbb G$ with $\mu = 0$
is exactly a one-dimensional Yetter-Drinfeld module over the 
group algebra $k[G]$, and the Hopf algebra $A(\mathbb G)$ is 
a Radford biproduct of the algebra $k[x]/(x^d)$ by $k[G]$.
More generally, this construction might be done for any one-dimensional
Yetter-Drinfeld module over an arbitrary Hopf algebra $H$:
see e.g. \cite{[AS],[CDMM]}.

\medskip

Let $\mathbb G = (G,g,\chi,\mu)$ be a group datum.
We need a few more informations concerning the Hopf
algebra $A(\mathbb G)$.
The grouplike elements  of $A(\mathbb G)$ coincide with the group $G$.
For $h \in G$, the set of 
$(1,h)$-primitives (i.e. elements $y \in A(\mathbb G)$
such that $\Delta(y) = 1 \otimes y + y \otimes h$)
is $kx \oplus k(g-1)$ if $h = g$ and $k(h-1)$ otherwise. 
From this we see that if $\mathbb G_1 = (G_1,g_1,\chi_1,\mu_1)$
and $\mathbb G_2 = (G_2,g_2,\chi_2,\mu_2)$ are group data, then
the Hopf algebras $A(\mathbb G_1)$ and $A(\mathbb G_2)$ are isomorphic
if and only if the group data are isomorphic in the sense of 
\cite{[CHYZ]}, namely there exists a group isomorphism
$f : G_1 \longrightarrow G_2$ such that $f(g_1) = g_2$ and 
$\chi_2 \circ f = \chi_1$, and $\delta \in \dot k$ such that
$\mu_1 = \delta^d \mu_2$. 
We also see that if $\mathbb G = (G,g,\chi)$ is a  
group datum of type I to V, then
$${\rm Aut}_{\rm Hopf}(A(\mathbb G))  \cong
{\rm Aut}_{g, \chi}(G) \times \dot k,$$
where ${\rm Aut}_{g, \chi}(G) = \{ u \in {\rm Aut}(G) \ | \
u(g)=g , \ \chi \circ u = \chi \}$.
If $\mathbb G = (G,g,\chi,\mu)$ is a type VI group datum,
then
$${\rm Aut}_{\rm Hopf}(A(\mathbb G))  \cong
{\rm Aut}_{g, \chi}(G) \times \mu_d.$$

\subsection{Cohomological preliminaries}

In this subsection we introduce modified versions
of the second cohomology group of a group. 
We will only need elementary group cohomology.
In this setting, an appropriate reference is \cite{[Kar]}.

Let $G$ be a group and let $g \in G$ be a central element. We put
$$Z^2_g(G,\dot k) = \{\sigma \in Z^2(G, \dot k), \ 
\sigma(g,h) = \sigma(h,g), \forall h \in G \} \quad {\rm and}$$  
$$B^2_g(G , \dot k) = \{ \partial(\mu), \
\mu : G \rightarrow \dot k, \ \mu(g)= 1 = \mu(1)\}.$$
For $g_1,g_2 \in Z(G)$, it is clear that
$B_{g_2}^2(G, \dot k)$ is a subgroup of $Z^2_{g_1}(G, \dot k)$ and
we define
$$H_{g_1,g_2}^2(G, \dot k) = Z^2_{g_1}(G,\dot k) / B_{g_2}^2(G,\dot k).$$
We have $ H_{1,1}^2(G, \dot k) =H^2(G, \dot k)$. 
These modified second cohomology groups will 
be useful for the description of Galois objects and biGalois 
groups over the non semisimple monomial Hopf algebras.
Some explicit computations will be done in Section 4.

Here is an another useful construction. Let $G$ be a group 
and let $g \in Z(G)$ be an element of order $n$.
Consider the group morphism
$$
\epsilon_0 :  Z^2(G,\dot k) \longrightarrow \dot k, \quad
\sigma \longmapsto \sigma(g,g) \ldots \sigma(g,g^{n-1}).$$ 
It is clear that $\epsilon_0$ induces a group morphism
$$\epsilon : H^2_{1,g}(G,\dot k) \longrightarrow \dot k.$$
In this way we have an action by automorphism 
of $H^2_{1,g}(G,\dot k)$ on the
additive group $(k,+)$, and we can form
the semi-direct product
$H^2_{1,g}(G,\dot k) \ltimes k$.

\medskip

Let us close this section by a useful fact, to be used freely in the
rest of the paper. Let $G$ be a group, let $H \subset G$
be a central subgroup and let $\sigma \in Z^2(G, \dot k)$.
Then the map
$$B_{\sigma} : G \times H \longrightarrow \dot k, \ 
(g,h) \longmapsto  \sigma(g,h)\sigma(h,g)^{-1},$$
is a pairing between the groups $H$ and $G$
(see \cite{[Kar]}, Lemma 2.2.6).

\section{Description of Galois objects}

In this section we determine the Galois objects for
the Hopf algebra associated to a group datum.
We only consider group data of type I,II, III, IV or V.
The type VI case will be treated thanks to the results in the 
next section, and may be reduced to the type III case: see Corollary 3.18.
Let us state the main result, for which the cases of type I, 
II, III, IV and V  will be proved in this section.

\begin{theo}
Let $\mathbb G = (G,g,\chi, \mu)$ be a group datum.
Then according to the type of $\mathbb G$,
we have the following description for ${\rm Gal}(A(\mathbb G))$.

\noindent
$\bullet$ Type I: ${\rm Gal}(A(\mathbb G)) \cong H^2(G,\dot k) \times k$.

\noindent
$\bullet$ Type II and IV: ${\rm Gal}(A(\mathbb G)) \cong H^2(G,\dot k)$.

\noindent
$\bullet$ Type III, V and VI: ${\rm Gal}(A(\mathbb G)) \cong H^2(G,\dot k) 
\amalg H^2_{g^d,g^d}(G, \dot k)$.
\end{theo}

The explicit bijections will be constructed during the proof.
Most of the proof will be done without any assumption 
on the type of the group datum. Case by case arguments
will be used only at the very end.

\medskip

We fix a group datum $\mathbb G = (G,g,\chi)$, with 
as usual $d = o(\chi(g)) >1$. 
The key construction is the following one.

\begin{defi}
For $\sigma \in Z^2(G, \dot k)$ and $a \in k$, we
define the algebra $A_{\sigma,a} (\mathbb G)$ 
to be  the algebra presented by generators 
$X,(T_h)_{h \in G}$ with defining relations,
$\forall h, h_1, h_2 \in G$:
$$T_{h_1}T_{h_2} = \sigma(h_1,h_2) T_{h_1h_2}, 
\quad T_1=1, \quad  XT_h = \chi(h) T_hX,
\quad X^d = a T_{g^d}.$$ 
\end{defi}

\begin{prop}
The algebra $A_{\sigma, a}(\mathbb G)$ has a right 
$A(\mathbb G)$-comodule algebra structure with coaction
$\alpha : A_{\sigma, a}(\mathbb G) \longrightarrow 
A_{\sigma, a}(\mathbb G) \otimes A(\mathbb G)$ defined by
$$\alpha(X) = 1 \otimes x + X \otimes g,
\quad \alpha(T_h) = T_h \otimes h, \quad \forall h \in G.$$
Then $A_{\sigma, a}(\mathbb G)$ is a right
$A(\mathbb G)$-Galois object if and only if
$$a \sigma(g^d,h) = a\chi(h)^d \sigma(h,g^d), \quad \forall h \in G.
\leqno(2.1)$$
In this case the set $\{T_h X^i, h \in G, 0 \leq i \leq d-1 \}$
is a linear basis of $A_{\sigma, a}(\mathbb G)$, and the
map $\Phi_{\sigma, a} : A(\mathbb G) 
\longrightarrow A_{\sigma, a}(\mathbb G)$, $hx^i \longmapsto T_h X^i$ 
is a right $A(\mathbb G)$-colinear isomorphism.
\end{prop}  

\noindent
\textbf{Proof}. It is easy to check that the map $\alpha$ 
is a well-defined algebra map that endows $A_{\sigma, a}(\mathbb G)$
with a right $A(\mathbb G)$-comodule algebra structure. 
Assume that $A_{\sigma, a}(\mathbb G)$ is a right 
$A(\mathbb G)$-Galois object. Then since $A(\mathbb G)$ is
finite-dimensional, it must be cleft and 
$\dim A_{\sigma, a}(\mathbb G) = |G|d$. The set 
$\{T_h X^i, h \in G, 0 \leq i \leq d-1 \}$ clearly
generate $A_{\sigma, a}(\mathbb G)$ linearly, and hence must
be a basis. Let $h \in G$. Then
$a \sigma (g^d,h) T_{g^dh} = a T_{g^d} T_h = X^d T_h
= \chi(h)^d T_h X^d = \chi(h)^d a T_h T_{g^d} =
\chi(h)^d a \sigma(h,g^d) T_{hg^d}$ and hence it follows
that (2.1) holds.

Conversely assume that (2.1) holds. 
Endow the set $S = \{X, T_h, h \in G \}$ with a total order such that
$T_h < X$, $\forall h \in G$. 
Then endow the set of monomials in elements of $S$ 
in such a way that two monomial of different
length are ordered according to their length, and that two
monomial of equal length are ordered lexicographically 
with respect to their indices. 
Our presentation is compatible with this order, there are three inclusion
ambiguities $(T_1T_h,T_1T_h)$, $(T_hT_1,T_hT_1)$ and
$(XT_1,XT_1)$ which are clearly resolvable, 
and there are the following overlap ambiguities:
$(T_{h_1}T_1,T_1)$,  $(T_1, T_1T_{h_2})$,  
$(T_{h_1}T_{h_2}, T_{h_2} T_{h_1})$,  
$(XT_h,T_hT_{h_2})$, $(XT_1,T_1)$, $(X^{i}X^{d-i}, X^{d-i}X^i)$, 
$(X^d,XT_h)$.
These ambiguities are easily seen to be resolvable,  
using (2.1) for the last one. We can use the diamond lemma
\cite{[Be]} to conclude that
$\{T_h X^i, h \in G, 0 \leq i \leq d-1 \}$
is a linear basis of $A_{\sigma, a}(\mathbb G)$.
It is clear that the map
$\Phi_{\sigma, a}$ 
is a right $A(\mathbb G)$-colinear isomorphism.
For $h \in G$, we have 
$$\kappa_r(\sigma(h^{-1},h)^{-1} T_{h^{-1}} \otimes T_h) = 1 \otimes h
\quad {\rm and}  \quad
\kappa_r(1 \otimes X - \sigma(g^{-1},g)^{-1} XT_{g^{-1}} \otimes T_g)
= 1 \otimes x$$
and since $x$ and the elements $h \in G$ generate
$A(\mathbb G)$ as an algebra, it is easy to see that $\kappa_r$ is
surjective. Hence $\kappa_r$ is an isomorphism and 
$A_{\sigma, a}(\mathbb G)$ is $A(\mathbb G)$-Galois. $\square$

\medskip

We want to prove that any right $A(\mathbb G)$-Galois object 
arises from the construction of Proposition 2.3. For this
we need a slightly more general construction.
Let $\sigma \in Z^2(G, \dot k)$, let $a \in k$
and let $\psi : G \longrightarrow k$.
We define the algebra $A_{\sigma,a, \psi} (\mathbb G)$ 
to be  the algebra presented by generators 
$X,(T_h)_{h \in G}$ with defining relations,
$\forall h, h_1, h_2 \in G$:
$$T_{h_1}T_{h_2} = \sigma(h_1,h_2) T_{h_1h_2}, 
\quad T_1=1, \quad  XT_h = \chi(h) T_hX + \psi(h) T_{gh},
\quad X^d = a T_{g^d}.$$ 

\begin{lemm}
The algebra $A_{\sigma, a, \psi}(\mathbb G)$ has a right 
$A(\mathbb G)$-comodule algebra structure with coaction
$\alpha : A_{\sigma, a,\psi}(\mathbb G) \longrightarrow 
A_{\sigma, a, \psi}(\mathbb G) \otimes A(\mathbb G)$ defined by
$$\alpha(X) = 1 \otimes x + X \otimes g,
\quad \alpha(T_h) = T_h \otimes h, \quad \forall h \in G.$$
Moreover the following assertions are equivalent.

\noindent
1)  $A_{\sigma,a, \psi} (\mathbb G)$ is a right 
$A(\mathbb G)$-Galois object.

\noindent
2) $\dim  A_{\sigma,a, \psi} (\mathbb G) = |G|d$ and the set
$\{T_hX^i, h \in G, 0 \leq i \leq d-1\}$ is a linear basis.

\noindent
3) $ A_{\sigma,a, \psi} (\mathbb G)$ is a non-zero algebra.

\end{lemm}  

\noindent
\textbf{Proof}. It is easy to check that the map $\alpha$ 
is a well-defined algebra map that endows $A_{\sigma, a, \psi}(\mathbb G)$
with a right $A(\mathbb G)$-comodule algebra structure. 
The proof of $1) \Rightarrow 2)$ is similar to the one in Proposition 2.3
and $2) \Rightarrow 3)$ is obvious. Similarly
to Proposition 2.3 the map $\kappa_r$ is surjective and if 
$A_{\sigma,a, \psi} (\mathbb G)$ is a non-zero algebra 
we have $\dim  A_{\sigma,a, \psi} (\mathbb G) \geq |G|d$,
hence $\dim  A_{\sigma,a, \psi} (\mathbb G) = |G|d$ and 
we conclude that $\kappa_r$ is an isomorphism. $\square$ 

\medskip

The following result motivates the introduction of the
algebras $ A_{\sigma,a, \psi} (\mathbb G)$.

\begin{lemm}
Let $T$ be a right $A(\mathbb G)$-Galois object.
Then there exists $\sigma \in Z^2(G, \dot k)$, $a \in k$
and $\psi : G \longrightarrow k$ such that
$T \cong A_{\sigma,a, \psi} (\mathbb G)$ as right
$A(\mathbb G)$-comodule algebras.
\end{lemm} 

\noindent
\textbf{Proof}. Let $\sigma : A(\mathbb G) \otimes  A(\mathbb G)
\longrightarrow k$ be a 2-cocycle such that
there exists an isomorphism of $ A(\mathbb G)$-comodule algebras
$f : {_{\sigma} \!A(\mathbb G)} \longrightarrow T$. We put 
$X = f(x)$ and $T_h = f(h)$ for $h \in G$. 
It is clear that for $h_1,h_2 \in G$, we have 
$T_{h_1}T_{h_2} = \sigma(h_1,h_2) T_{h_1h_2}$.
Let $\beta$ be the coaction of $A(\mathbb G)$ on $T$.
We have 
$$\beta(T_h) = T_h \otimes h , \ \forall h \in G \quad 
 {\rm and} \quad  \beta(X) = 1 \otimes x + X \otimes g.$$
Hence for $h \in G$, we have
$$\beta(XT_h - \chi(h)T_hX) = (XT_h - \chi(h)T_hX) \otimes hg,$$ 
and thus there exists $\psi(h) \in k$ such that
$XT_h - \chi(h) T_hX = \psi(h) T_{gh}$. This defines
a map $\psi : G \rightarrow k$. Similarly
$\beta(X^d) = (1 \otimes x + X\otimes g)^d = X^d \otimes g^d$ and
hence there exists $a \in k$ such that $X^d = aT_{g^d}$.
In this way we have a surjective right $A(\mathbb G)$-comodule
algebras morphism $A_{\sigma, a, \psi}(\mathbb G) \longrightarrow T$ 
where  $\sigma$ is the restriction of $\sigma$ to $G$.
Hence $ A_{\sigma, a ,\psi}(\mathbb G)$ is a non-zero algebra,
is  $A(\mathbb G)$-Galois by the previous lemma, 
and since $A_{\sigma,a ,\psi}(\mathbb G)$ and $T$ 
are both $A(\mathbb G)$-Galois, they are isomorphic via $f$. $\square$

\medskip

Therefore we need to study when the algebra $A_{\sigma, a, \psi}(\mathbb G)$
is $A(\mathbb G)$-Galois. The technical conditions are given in the
following result.

\begin{lemm}
If $A_{\sigma, a, \psi}(\mathbb G)$ is a right
$A(\mathbb G)$-Galois object, then for $h_1, h_2, h  \in G$
we have
$$\psi(h_1h_2) = 
\sigma(h_1,h_2)^{-1}( \chi(h_1) \sigma(h_1,h_2g) \psi(h_2)
+ \sigma(h_1g,h_2) \psi(h_1)),
\leqno(2.2)$$
$$\psi(h) = \psi(g) (\sigma(g,h) - \chi(h)\sigma(h,g))
\sigma(g,g)^{-1} (1 - \chi(g))^{-1}, \leqno(2.3)$$
$$a(\sigma(g^d,h)-\chi(h)^d\sigma(h,g^d)) = \psi(h) \psi(gh)
\ldots \psi(g^{d-1}h). \leqno(2.4)$$
\end{lemm}

\noindent
\textbf{Proof}.
Let $h_1,h_2 \in G$. We have
$$XT_{h_1h_2} = \chi(h_1h_2) T_{h_1h_2} X + \psi(h_1h_2)T_{h_1h_2g}.$$
On the other hand we have
\begin{align*}
X& T_{h_1h_2}  = \sigma(h_1,h_2)^{-1} XT_{h_1}T_{h_2} =
\sigma(h_1,h_2)^{-1}(\chi(h_1)T_{h_1}X + \psi(h_1) T_{h_1g}) T_{h_2} = \\
= & \sigma(h_1,h_2)^{-1}(\chi(h_1)\chi(h_2)T_{h_1}T_{h_2}X
+ \chi(h_1) \psi(h_2)T_{h_1}T_{h_2g} + \psi(h_1) T_{h_1g}T_{h_2}) = \\
= &  \chi(h_1h_2) T_{h_1h_2}X + 
\sigma(h_1,h_2)^{-1}(\chi(h_1)\psi(h_2) \sigma(h_1,h_2g)
+ \psi(h_1) \sigma(h_1g,h_2)) T_{h_1h_2g}.
\end{align*}
Hence by Lemma 2.3, Equation (2.2) holds.  

Let $h \in G$. Since $g$ is central in $G$, we have 
$\psi(hg) = \psi(gh)$, and it is a straightforward computation
to check that Equation (2.3) holds,
using (2.2) and the fact that $\sigma$ is a 2-cocycle.

For $n  \in \mathbb N^*$, there exists $a_1, \ldots a_{n-1} \in k$
such that
$$X^n T_h = \chi(h)^n T_h X^n +
\sum_{l = 1}^{n-1} a_l T_{g^lh} X^{n-l} + \psi(h) \psi(gh) \ldots
\psi(g^{n-1}h)T_{g^nh}.$$
This is proved easily by induction.
For $n=d$, we see, using Lemma 2.3, that Equation (2.4) holds. 
$\square$

\begin{lemm}
Let $\sigma \in Z^2(G,\dot k)$. Then for $h \in G$, we have
$$\prod_{i=0}^{d-1}
\left(\sigma(g,g^ih) - \chi(g^ih)\sigma(g^ih,g)\right) =
\sigma(g,g) \ldots \sigma(g,g^{d-1}) \left(\sigma(g^d,h)
- \chi(h)^d\sigma(h,g^d)\right).$$
\end{lemm}

\noindent
\textbf{Proof}. Consider the polynomials 
$$P(X) = \prod_{i=0}^{d-1}
\left(\sigma(g,g^ih) - \chi(g^i)\sigma(g^ih,g)X\right) 
\in k[X] \quad {\rm and}$$
$$Q(X) = 
\sigma(g,g) \ldots \sigma(g,g^{d-1}) \left(\sigma(g^d,h)
- \sigma(h,g^d)X^d\right) \in k[X].$$
Since  $\sigma$ is a 2-cocycle and since $g$ is central, 
we have  for $i \in \mathbb N$:
$$\sigma(g,h)\sigma(h,g)^{-1} = \sigma(g,g^ih)\sigma(g^ih,g)^{-1}
\ {\rm and} \
(\sigma(g,h)\sigma(h,g)^{-1})^i = \sigma(g^i,h)\sigma(h,g^i)^{-1}.$$
Let $t$ be a root of $P$: this means that there exists
$i \in \{0, \ldots , d-1\}$ such that 
$t = \sigma(g,g^ih) \sigma(g^ih,g)^{-1} \chi(g)^{-i}=
\sigma(g,h)\sigma(h,g)^{-1}\chi(g)^{-i}$.
Then 
$t^d = (\sigma(g,h) \sigma(h,g)^{-1})^d =
\sigma(g^d,h) \sigma(h,g^d)^{-1}$.
Hence $t$ is a root of $Q$.
Conversely assume that $t$ is a root of $Q$. Then
we have $t^d \sigma(g^d,h)^{-1} \sigma(h,g^d)  = 1
= (t \sigma(g,h)^{-1} \sigma(h,g))^d$.
Hence there exists $i \in \{0, \ldots , d-1 \}$ such that
$t = \chi(g)^{-i} \sigma(g,h) \sigma(h,g)^{-1} =
\chi(g)^{-i} \sigma(g,g^ih) \sigma(g^ih,g)^{-1}$ and 
$t$ is a root of $P$.
It is easy to see that $P$ and $Q$ have the same constant term,
and we conclude that $P=Q$. In particular $P(\chi(h)) = Q(\chi(h))$:
this proves the lemma. $\square$

\begin{lemm}
Let $\sigma \in Z^2(G,\dot k)$, $a \in k$ and
$\psi : G \rightarrow k$ be such that 
$A_{\sigma, a,\psi}(\mathbb G)$ is $A(\mathbb G)$-Galois.
Then there exists $a' \in k$ such that
$A_{\sigma, a,\psi}(\mathbb G) \cong A_{\sigma, a'}(\mathbb G)$
as right $A(\mathbb G)$-comodule algebras.
\end{lemm}

\noindent
\textbf{Proof}. We put
$$a' = a - \psi(g)^d \sigma(g,g)^{-d} (1-\chi(g))^{-d}
\sigma(g,g) \dots \sigma(g,g^{d-1}).$$
Let us first show that
$A_{\sigma,a'}$ is $A(\mathbb G)$-Galois. Let $h \in G$.
We have
\begin{align*}
a' & \left(\sigma(g^d,h) - \chi(h)^d\sigma(h,g^d)\right) = 
a\left(\sigma(g^d,h) - \chi(h)^d\sigma(h,g^d)\right) \\
- & \psi(g)^d \sigma(g,g)^{-d} (1-\chi(g))^{-d}
\sigma(g,g) \dots \sigma(g,g^{d-1}) 
\left(\sigma(g^d,h) - \chi(h)^d\sigma(h,g^d)\right) \\
= & \psi(h) \psi(gh) \ldots \psi(g^{d-1}h) \ 
({\rm by } \ {\rm  Lemma} \ 2.6) \  - \\ 
\psi(g)^d & \sigma(g,g)^{-d} (1-\chi(g))^{-d}
\sigma(g,g) \dots \sigma(g,g^{d-1}) 
\left(\sigma(g^d,h) - \chi(h)^d\sigma(h,g^d)\right) \\
= &  \psi(g)^d  \sigma(g,g)^{-d} (1-\chi(g))^{-d}  \\
[(\prod_{i=0}^{d-1} &
\left(\sigma(g,g^ih) - \chi(g^ih)\sigma(g^ih,g)\right) -
\sigma(g,g) \ldots \sigma(g,g^{d-1}) \left(\sigma(g^d,h)
- \chi(h)^d\sigma(h,g^d)\right)] \\
(& {\rm by} \ {\rm  Lemma} \ 2.6) \ = 0 \ ({\rm by} \ {\rm  Lemma} \ 2.7).
\end{align*}
Hence $A_{\sigma,a'}(\mathbb G)$ is $A(\mathbb G)$-Galois
by Proposition 2.3. 
Now put $\lambda = \psi(g) \sigma(g,g)^{-1} (1-\chi(g))^{-1}$. 
It is straightforward to check, using Lemma 2.6, that
there exists an $A(\mathbb G)$-comodule algebra morphism 
$f : A_{\sigma,a, \psi}(\mathbb G)
\longrightarrow A_{\sigma,a'}(\mathbb G)$ such that
$f(X) = X + \lambda T_g$ and $f(T_h) = T_h$, $\forall h \in G$.
Then $f$ is an isomorphism 
since $A_{\sigma,a, \psi}(\mathbb G)$ and $A_{\sigma,a'}(\mathbb G)$
are $A(\mathbb G)$-Galois. $\square$ 

\medskip

Combining Lemmas 2.5 and 2.8, we get the following result.

\begin{prop}
Let $T$ be a right $A(\mathbb G)$-Galois object.
Then there exists $\sigma \in Z^2(G, \dot k)$ and  $a \in k$
such that
$T \cong A_{\sigma,a} (\mathbb G)$ as right $A(\mathbb G)$-comodule
algebras.
\end{prop}

The next step is to classify the Galois objects $A_{\sigma, a}(\mathbb G)$. 

\begin{prop}
Let $\sigma, \tau \in Z^2(G,\dot k)$ and $a,b \in k$
be such that $A_{\sigma, a}(\mathbb G)$ and $A_{\tau, b}(\mathbb G)$
are $A(\mathbb G)$-Galois objects. Then the right
$A(\mathbb G)$-comodule algebras $A_{\sigma, a}(\mathbb G)$ and
$A_{\tau, b}(\mathbb G)$ are isomorphic if and only if
there exists $\mu : G \longrightarrow \dot k$ with $\mu(1) =1$ such
that
$$\sigma = \partial(\mu) \tau \quad {\rm and} \quad 
b = a \mu(g^d).$$ 
\end{prop}

\noindent
\textbf{Proof}. Let $f : A_{\sigma, a}(\mathbb G)
\longrightarrow A_{\tau, b}(\mathbb G)$ be an $A(\mathbb G)$-colinear algebra
isomorphism. Then $\Phi_{\tau,b}^{-1} \circ f \circ \Phi_{\sigma, a}$
is a right $A(\mathbb G)$-colinear automorphism of $A(\mathbb G)$, and
hence there exists a convolution invertible linear map
$\phi : A(\mathbb G) \rightarrow k$ such that 
$\Phi_{\tau,b}^{-1} \circ f \circ \Phi_{\sigma, a} = \phi * {\rm id}$.
From this we see easily see that there exists
$\mu : G \longrightarrow \dot k$ and $\lambda \in k$ such that
$$f(X) = X + \lambda T_g \quad {\rm and} \quad
f(T_h) = \mu(h) T_h, \ \forall h \in G.$$ 
It is then clear that $\sigma = \partial(\mu) \tau$ and that
$\mu(1) = 1$.
We have
$$f(XT_g) = \mu(g) XT_g + \lambda \mu(g) \tau(g,g) T_{g^2}
\quad {\rm and}$$
$$f(\chi(g) T_g X)= \chi(g) \mu(g) T_g X
+ \chi(g) \lambda \tau(g,g) \mu(g)T_{g^2}.$$
Since $\chi(g) \not = 1$, we have $\lambda = 0$.
Then $bT_{g^d} = X^d = f(X^d) = \mu(g^d) a T_{g^d}$ and hence
$b = a \mu(g^d)$. 

Conversely if  $\sigma = \partial(\mu) \tau$ 
with $\mu(1) = 1$ and $b = a \mu(g^d)$, 
it is straightforward to check that there exists
a right $A(\mathbb G)$-comodule algebra isomorphism 
$f :  A_{\sigma, a}(\mathbb G)
\rightarrow A_{\tau, b}(\mathbb G)$ such that
$f(X) = X$ and $f(T_h) = \mu(h)T_h$, $\forall h \in G$. $\square$

\medskip

We can now finish the proof of Theorem 2.1 with case by case
arguments.

\begin{prop}
Let $\mathbb G = (G,g, \chi)$ be a type I group
datum. For any $\sigma \in Z^2(G, \dot k)$ and $a \in k$, 
then $A_{\sigma, a}(\mathbb G)$ is $A(\mathbb G)$-Galois.
The map $Z^2(G,\dot k) \times k \longrightarrow {\rm Gal}(A(\mathbb G))$,
$(\sigma,a) \longmapsto [A_{\sigma,a}(\mathbb G)]$, 
induces a bijection
$$H^2(G, \dot k) \times k \cong {\rm Gal}(A(\mathbb G)).$$
\end{prop} 

\noindent
\textbf{Proof}. It is clear that for any $\sigma \in Z^2(G, \dot k)$ 
and $a \in k$, condition (2.1) is satisfied and hence 
$A_{\sigma, a}(\mathbb G)$ is $A(\mathbb G)$-Galois by 
Proposition 2.3. The map 
$Z^2(G,\dot k) \times k \longrightarrow {\rm Gal}(A(\mathbb G))$,
$(\sigma,a) \longmapsto [A_{\sigma,a}(\mathbb G)]$, 
induces, by Proposition 2.10, an injective map
$H^2(G, \dot k) \times k \longrightarrow {\rm Gal}(A(\mathbb G))$.
This map is surjective by Proposition 2.9: this concludes the proof.
$\square$

\begin{prop}
Let $\mathbb G = (G,g, \chi)$ be a type II group
datum. For any $\sigma \in Z^2(G, \dot k)$, 
then $A_{\sigma, 0}(\mathbb G)$ is $A(\mathbb G)$-Galois.
The map $Z^2(G,\dot k) \longrightarrow {\rm Gal}(A(\mathbb G))$,
$\sigma \longmapsto [A_{\sigma,0}(\mathbb G)]$, 
induces a bijection
$$H^2(G, \dot k) \cong {\rm Gal}(A(\mathbb G)).$$
\end{prop} 

\noindent
\textbf{Proof}. For any $\sigma \in Z^2(G, \dot k)$ 
and $a =0$, condition (2.1) is satisfied and hence 
$A_{\sigma, 0}(\mathbb G)$ is $A(\mathbb G)$-Galois by 
Proposition 2.3. The map 
$Z^2(G,\dot k) \longrightarrow {\rm Gal}(A(\mathbb G))$,
$\sigma \longmapsto [A_{\sigma,0}(\mathbb G)]$, 
induces, by Proposition 2.10, an injective map
$H^2(G, \dot k) \longrightarrow {\rm Gal}(A(\mathbb G))$.
Let $T$ be a right $A(\mathbb G)$-Galois object: by Proposition 2.10
there exists $\sigma \in Z^2(G, \dot k)$ and $a \in k$
such that $T \cong A_{\sigma,a}(\mathbb G)$. But since 
$\chi^d \not = 1$, condition (2.1) is satisfied only if $a=0$.
This concludes the proof. $\square$

\begin{prop}
Let $\mathbb G = (G,g, \chi)$ be a type III group
datum. Then
the maps $Z^2(G,\dot k) \longrightarrow {\rm Gal}(A(\mathbb G))$,
$\sigma \longmapsto [A_{\sigma,0}(\mathbb G)]$,
and   $Z^2_{g^d}(G,\dot k) \longrightarrow {\rm Gal}(A(\mathbb G))$, \break
$\sigma \longmapsto [A_{\sigma,1}(\mathbb G)]$,
induce a bijection
$$H^2(G, \dot k) \amalg
 H^2_{g^d,g^d}(G,\dot k) \cong {\rm Gal}(A(\mathbb G)).$$
\end{prop} 

\noindent
\textbf{Proof}. Our two maps induce a map
$H^2(G, \dot k) \amalg H^2_{g^d,g^d}(G,\dot k) \longrightarrow
 {\rm Gal}(A(\mathbb G))$, which is injective by Proposition 2.10.
Let $T$ be a right $A(\mathbb G)$-Galois object: by Proposition 2.9
there exists $\sigma \in Z^2(G, \dot k)$ and $a \in k$
such that $T \cong A_{\sigma,a}(\mathbb G)$.
If $a \not = 0$, then $\sigma \in Z^2_{g^d}(G,\dot k)$ 
by Proposition 2.3. 
Consider $\mu : G \longrightarrow \dot k$ such that
$\mu(1)=1$ and $\mu(g^d)= a$. Then 
$A_{\sigma,a}(\mathbb G) \cong 
A_{\partial(\mu)\sigma,1}(\mathbb G)$ by Proposition 2.10, and 
hence our map is also surjective. $\square$

\medskip

The proofs of the type IV and V cases are similar and left to the reader,
and conclude the proof of Theorem 2.1 (up to the type VI case).

\begin{prop}
Let $\mathbb G = (G,g, \chi)$ be a type IV group
datum. Then 
the map $Z^2(G,\dot k) \longrightarrow {\rm Gal}(A(\mathbb G))$,
$\sigma \longmapsto [A_{\sigma,0}(\mathbb G)]$,
induces a bijection
$$H^2(G, \dot k)  
\cong {\rm Gal}(A(\mathbb G)).$$
\end{prop} 

\begin{prop}
Let $\mathbb G = (G,g, \chi)$ be a type V group datum. 
Let $\tau \in Z^2(G,\dot k)$ be such that
$\chi(h)^d = \tau(g^d,h)\tau(h,g^d)^{-1}$, $\forall h \in G$.
Then $\forall \sigma \in Z^2_{g^d}(G,\dot k)$, 
$A_{\tau\sigma,1}(\mathbb G)$ is $A(\mathbb G)$-Galois. 
The maps $Z^2(G,\dot k) \longrightarrow {\rm Gal}(A(\mathbb G))$,
$\sigma \longmapsto [A_{\sigma,0}(\mathbb G)]$,
and   \break $Z^2_{g^d}(G,\dot k) \longrightarrow {\rm Gal}(A(\mathbb G))$,
$\sigma \longmapsto [A_{\tau\sigma,1}(\mathbb G)]$,
induce a bijection
$$H^2(G, \dot k) \amalg H^2_{g^d,g^d}(G, \dot k) 
\cong {\rm Gal}(A(\mathbb G)).$$
\end{prop} 

We conclude the section with  a few words concerning
the notion of homotopy equivalence of Hopf-Galois extensions,
recently introduced by Kassel and Schneider \cite{[KS]}.
Let $\mathbb G = (G,g,\chi,\mu)$ be a group datum
and let $\mathcal H_k(A(\mathbb G))$ be the set of homotopy
classes of right $A(\mathbb G)$-Galois extensions of $k$
(see \cite{[KS]} for the precise definition). Then we have
$$\mathcal H_k(A(\mathbb G)) \cong H^2(G,\dot k).$$
Indeed if $\mu = 0$, then 
$A(\mathbb G) = \oplus_{i=0}^{d-1} k[G]x^i$ and hence $A(\mathbb G)$ is
an $\mathbb N$-graded Hopf algebra. Therefore
combining Corollary 2.9 and Proposition 3.2 of \cite{[KS]},
we have the claimed result. Then one gets the result in the
type VI case by combining  Corollary 2.11 of \cite{[KS]}
and Corollary 3.18 in the next section.

\section{BiGalois groups}

In this section we give a general description
of the biGalois group of a monomial Hopf algebra.
Several particular cases will be studied
in more detail in the next section.

Just like in the previous section, our first goal is to state our
main result. We first need to introduce
some terminology. Let 
$\mathbb G = (G,g, \chi, \mu)$ be a group datum. 
We consider the subgroup Aut$_g(G)$ of elements  
$u \in {\rm Aut}(G)$ satisfying $u(g) = g$ and 
the modified cohomology group 
$H^2_{1,g}(G, \dot k) = Z^2(G, \dot k)/B_g^2(G, \dot k)$
of the first section. The group Aut$_g(G)$ has a natural
right action by automorphisms on $H^2_{1,g}(G, \dot k)$, and hence
we may form the semi-direct product
${\rm Aut}_g(G) \ltimes H^2_{1,g}(G, \dot k)$. We define
now
$$\Gamma(\mathbb G)= \{
(u,\overline{\sigma}) \in {\rm Aut}_g(G) \times H^2_{1,g}(G, \dot k), \ 
\chi \circ u(h) = \sigma(g,h)^{-1} \sigma(h,g) \chi(h), \ 
\forall h \in G \}.$$    
It is clear that, since $g$ is central in $G$, that the defining
equation for the elements of $\Gamma(\mathbb G)$ does not depend
of the choice of a representant of a (modified) cohomology class.
The following lemma is a straightforward verification.

\begin{lemm}
Let $\mathbb G = (G,g, \chi, \mu)$ be a  group datum. Then
$\Gamma(\mathbb G)$ is a subgroup of \break
${\rm Aut}_g(G) \ltimes H^2_{1,g}(G, \dot k)$.
\end{lemm}

The group morphism 
$\epsilon : H^2_{1,g}(G, \dot k) \longrightarrow \dot k$
has a natural continuation to 
$\Gamma(\mathbb G)$, and hence we may form the semi-direct 
product $\Gamma (\mathbb G) \ltimes k$. We have now all the
ingredients to state the main result of the section.

\begin{theo}
Let $\mathbb G = (G,g, \chi)$ be a group datum. Then according
to the type of $\mathbb G$, we have the following group isomorphism.

\noindent 
$\bullet$ Type I : ${\rm BiGal}(A(\mathbb G)) \cong
\Gamma(\mathbb G) \ltimes k$.

\noindent 
$\bullet$ Type II, III, IV, V and VI: 
${\rm BiGal}(A(\mathbb G)) \cong \Gamma(\mathbb G)$.
\end{theo}
 
Types V and VI will be shown to reduce to type III.
Similarly to the previous section, we begin
with general results and constructions. 

\begin{prop} 
Let $\mathbb G = (G,g, \chi)$ be a group datum.
Let $\sigma \in Z^2(G, \dot k)$ and $a \in k$ 
with $a=0$ if $\mathbb G$ is not of type I: 
$A_{\sigma,a}(\mathbb G)$ is $A(\mathbb G)$-Galois.
Assume that there exists $u \in {\rm Aut}_g(G)$
such that $(u, \overline{\sigma}) \in \Gamma(\mathbb G)$.
Then $A_{\sigma,a}(\mathbb G)$ has a left $A(\mathbb G)$-comodule
algebra structure $\beta_u : 
A_{\sigma,a}(\mathbb G)\longrightarrow A(\mathbb G)
\otimes A_{\sigma,a}(\mathbb G)$ defined by
$$\beta_u(X) = 1 \otimes X + x \otimes T_g, \quad 
\beta_u(T_h) = u(h) \otimes T_h, \ \forall h \in G.$$
Furthermore $A_{\sigma,a}(\mathbb G)$ is an 
$A(\mathbb G)$-$A(\mathbb G)$-bicomodule algebra and 
is $A(\mathbb G)$-biGalois. 
This Hopf biGalois extension is denoted by 
$A_{\sigma,a}^u(\mathbb G)$.
\end{prop}

\noindent
\textbf{Proof}. It is easy to check that $\beta_u$
is a well-defined algebra map, using that 
$(u, \overline{\sigma}) \in \Gamma(\mathbb G)$ when
checking that $\beta_u(XT_h) = \beta_u(\chi(h) T_hX)$, 
that endows $A_{\sigma, a}(\mathbb G)$ with an 
$A(\mathbb G)$-$A(\mathbb G)$-bicomodule structure.
There remains to check that
$\kappa_l : A_{\sigma, a}(\mathbb G) \otimes A_{\sigma, a}(\mathbb G)
\longrightarrow A(\mathbb G) \otimes A_{\sigma, a}(\mathbb G)$ is bijective.
We have
$\kappa_l(\sigma(g,g^{-1})^{-1}(X \otimes T_{g^{-1}} -
1 \otimes XT_{g^{-1}})) = x \otimes 1$ and 
$\kappa_l(\sigma(h,h^{-1})^{-1} T_h \otimes T_{h^{-1}})
= u(h) \otimes 1$ for $h \in G$.
Since the elements $x$ and $u(h)$, $h \in G$ generate 
$A(\mathbb G)$ as an algebra and since $\kappa_l$ is right 
$A_{\sigma, a}(\mathbb G)$-linear, we conclude that it is surjective 
and an isomorphism. $\square$.

\medskip

The next step is to classify the Hopf biGalois
extensions just constructed.

\begin{prop} 
Let $\mathbb G = (G,g, \chi)$ be a group datum.
Let $\sigma , \tau \in Z^2(G, \dot k)$ and $a,b \in k$ 
with $a=b=0$ if $\mathbb G$ is not of type I: 
$A_{\sigma,a}(\mathbb G)$ and $A_{\tau,b}(\mathbb G)$ are  
$A(\mathbb G)$-Galois.
Assume that there exist $u,v \in {\rm Aut}_g(G)$
such that $(u, \overline{\sigma}) \in \Gamma(\mathbb G)$ and
$(v, \overline{\tau}) \in \Gamma(\mathbb G)$.
Then the $A(\mathbb G)$-biGalois extensions
$A_{\sigma,a}^u(\mathbb G)$ and $A_{\tau,b}^v(\mathbb G)$
are isomorphic if and only if $u=v$ and there exists
$\mu : G \longrightarrow \dot k$ satisfying
$\mu(1) = \mu(g) = 1$, $\sigma = \partial(\mu) \tau$
and  $b = a \mu(g^d)$.
\end{prop}

\noindent
\textbf{Proof}. Let $f : A_{\sigma, a}^u(\mathbb G) \longrightarrow
A_{\tau, b}^v(\mathbb G)$ be an $A(\mathbb G)$-bicolinear
isomorphism. By Proposition 2.10 and its proof, there
exists $\mu : G \longrightarrow \dot k$ with $\mu(g) = 1$,
$\sigma = \partial(\mu) \tau$, $b= a \mu(g^d)$, such that
$f(X) = X$ and $f(T_h) = \mu(h) T_h$ for $h \in G$.
Using that $f$ is left colinear, we see
that $u=v$ and that $ \mu(g) =1$.
The converse assertion is clear from the above considerations.
$\square$

\medskip

We now connect the group structures of 
$\Gamma(\mathbb G)$ and BiGal$(A(\mathbb G))$.

\begin{prop} 
Let $\mathbb G = (G,g, \chi)$ be a group datum.
Let $\sigma , \tau \in Z^2(G, \dot k)$ and $a,b \in k$ 
with $a=b=0$ if $\mathbb G$ is not of type I: 
$A_{\sigma,a}(\mathbb G)$ and $A_{\tau,b}(\mathbb G)$ 
are  $A(\mathbb G)$-Galois. 
Assume that there exist $u,v \in {\rm Aut}_g(G)$
such that $(u, \overline{\sigma}) \in \Gamma(\mathbb G)$
and $(v, \overline{\tau}) \in \Gamma(\mathbb G)$.
Put $\nu = (\sigma \circ v \times v)\tau$ and 
$c = \tau(g,g) \ldots \tau(g,g^{d-1})a+b$.
Then $A_{\nu,c}(\mathbb G)$ is right $A(\mathbb G)$-Galois,
$(u \circ v, \overline{\nu}) \in \Gamma(\mathbb G)$ and
the $A(\mathbb G)$-biGalois extensions 
$A_{\nu,c}^{u\circ v}(\mathbb G)$ and
$A_{\sigma,a}^u(\mathbb G) \square_{A(\mathbb G)} A_{\tau,b}^v(\mathbb G)$
are isomorphic. 
\end{prop}

\noindent
\textbf{Proof}. It is clear that
$A_{\nu,c}(\mathbb G)$ is $A(\mathbb G)$-Galois.
We have 
$(u \circ v, \overline{\nu}) \in \Gamma(\mathbb G)$ by Lemma 3.1 and 
hence we can consider the  $A(\mathbb G)$-biGalois extension
 $A_{\nu,c}^{u\circ v}(\mathbb G)$. It is easy to see that
$1 \otimes X + X \otimes T_g \in A_{\sigma,a}^u(\mathbb G) \square_{A(\mathbb G)} A_{\tau,b}^v(\mathbb G)$ and that for $h \in G$,
$T_{v(h)} \otimes T_h \in  A_{\sigma,a}^u(\mathbb G) \square_{A(\mathbb G)} A_{\tau,b}^v(\mathbb G)$. 
One checks now that there exists an algebra morphism
$$\gamma : A_{\nu, c}^{u \circ v}(\mathbb G) \longrightarrow
A_{\sigma,a}^u(\mathbb G) \square_{A(\mathbb G)} A_{\tau,b}^v(\mathbb G)$$
such that $\gamma (X) = 1 \otimes X + X \otimes T_g$ and
$\gamma(T_h) = T_{v(h)} \otimes T_h$, $\forall h \in G$.
It is clear that $\gamma$ is $A(\mathbb G)$-bicolinear, and since
it is a morphism of Galois extensions, it is an isomorphism. $\square$ 

\medskip

Now we need to show that we have constructed all 
the $A(\mathbb G)$-biGalois extensions.
We will need case by case arguments.
We begin with a property which, unfortunately,
only holds in general for type I or II group data.

\begin{lemm}
Let $\mathbb G = (G,g, \chi)$ be a group datum of type I or II.
Let $\sigma \in Z^2(G, \dot k)$. Then there exists
$u \in {\rm Aut}_g(G)$ such that $(u, \overline{\sigma}) \in 
\Gamma(\mathbb G)$.
\end{lemm}

\noindent
\textbf{Proof}. Since $\mathbb G$ is of type I or II we
have $d = o(\chi(g)) = o(g)$ and 
$\chi$ induces a group isomorphism $\langle g \rangle \cong \mu_d$.
Let $s : \mu_d \longrightarrow \langle g \rangle$ be 
a group morphism such that
$\chi \circ s = {\rm id}$. Now let $\psi : G \longrightarrow \dot k$
be a group morphism. Define a map $u : G \longrightarrow G$
by $u(h) = s(\psi(h)) h$, $\forall h \in G$. Then it is easy to
see that $u$ is a group morphism ($g$ is central) and
that $\chi \circ u (h) = \psi(h) \chi(h)$, $\forall h \in G$.
Furthermore $u$ is an isomorphism and $u(g)=g$ if $\psi(g)=1$.
This general argument, applied to the group morphism
$G \longrightarrow \mu_d$, $h \longmapsto 
\sigma(g,h)^{-1} \sigma(h,g)$, proves our lemma. $\square$ 

\medskip

Combining Propositions 2.11, 2.12 and  3.3, Lemma 3.6
with Theorem 3.5 and Corollary 5.7 of \cite{[Sc1]},
we get the following result, which generalizes
Corollaries 4 and 18 in \cite{[Sc2]}.

\begin{coro}
Let $\mathbb G = (G,g, \chi)$ be a group datum of type I or II.
Then any right $A(\mathbb G)$-Galois object is $A(\mathbb G)$-biGalois.
The Hopf algebra $A(\mathbb G)$ is categorically rigid: 
if $H$ is any Hopf algebra such that the monoidal $k$-linear
categories of comodules of $A(\mathbb G)$ and $H$ are equivalent, then
$H \cong A(\mathbb G)$ as Hopf algebras.
\end{coro}

We can finish the proof of Theorem 3.2 in the case of type 
I and II group data.

\begin{prop}
Let $\mathbb G = (G,g, \chi)$ be a type I group datum.
Then we have a group isomorphism:
\begin{align*}
\Psi : \Gamma(\mathbb G) \ltimes k 
\longrightarrow & {\rm BiGal}(A(\mathbb G)) \\
(u, \overline{\sigma},a ) \longmapsto & [A_{\sigma,a}^u].
\end{align*}
\end{prop}

\noindent
\textbf{Proof}. It follows immediately from 
Propositions 2.11, 3.3 and 3.4 that $\Psi$ is a well-defined injective 
map. Also it is clear from Proposition 3.5 that $\psi$ is
a group morphism. There just remains to show that
$\Psi$ is surjective. Let $Z$ be an $A(\mathbb G)$-biGalois extension.
By Proposition 2.11 we can assume that 
$Z = A_{\sigma,a}(\mathbb G)$ as right $A(\mathbb G)$-Galois object.
By Lemma 3.6 there exists $u \in {\rm Aut}_g(G)$ such that 
$(u, \overline{\sigma}) \in \Gamma(\mathbb G)$, and hence
by Proposition 3.3 $A_{\sigma,a}^u(\mathbb G)$ is $A(\mathbb G)$-biGalois.
By  \cite{[Sc1]}, Theorem 3.5, there exists 
$f \in {\rm Aut}_{\rm Hopf}(A(\mathbb G))$ such
that $Z = {^f} \!A_{\sigma,a}^u(\mathbb G)$ as $A(\mathbb G)$-bicomodule
algebras, i.e. if $\beta '$ denotes the left coaction
of $Z$, we have $\beta ' = (f \otimes {\rm id}) \circ \beta_u$. 
 There exists $v \in {\rm Aut}_{g, \chi}(G)$ and $r \in \dot k$
such that $f(h) = v(h)$, $\forall h \in G$ and $f(x) = rx$.
Consider now a map $\mu : G \longrightarrow \dot k$
with $\mu(1) =1$ and $\mu(g) = r$. It is clear that
$(v \circ u, \overline{\partial(\mu)\sigma}) \in \Gamma(\mathbb G)$, and
hence we consider the $A(\mathbb G)$-biGalois extension
$A_{\partial(\mu)\sigma,a}^{v \circ u}(\mathbb G)$.
It is not difficult to check that the right
$A_{\sigma,a}(\mathbb G)$-comodule algebra isomorphism
$f' : A_{\partial(\mu)\sigma,a}^{v \circ u}(\mathbb G) \longrightarrow
{^f} \! A_{\sigma,a}^u(\mathbb G)$
defined by $f'(X) = X$ and $f'(T_h) = \mu(h)T_h$, $\forall h \in G$,
is also left colinear. Therefore $\Psi$ is also surjective $\square$   

\medskip

The type I case will be examined in more detail in the next section.
The proof of the type II case is similar and left to the reader.

\begin{prop}
Let $\mathbb G = (G,g, \chi)$ be a type II group datum.
Then we have a group isomorphism:
\begin{align*}
\Psi : \Gamma(\mathbb G) 
\longrightarrow & {\rm BiGal}(A(\mathbb G)) \\
(u, \overline{\sigma}) \longmapsto & [A_{\sigma,0}^u].
\end{align*}
\end{prop}

We now begin the study of the biGalois group in  the
type IV case.

\begin{prop} 
Let $\mathbb G = (G,g, \chi)$ be a type IV group datum.
Let $\sigma \in Z^2(G, \dot k)$.
Consider the group morphism $\chi ' : G \longrightarrow \dot k$
defined by $\chi'(h) = \sigma(g,h)^{-1} \sigma(h,g) \chi(h)$.
Then the group datum $\mathbb G' = (G,g , \chi ')$ is a type IV group
datum. The right $A(\mathbb G)$-Galois object
$A_{\sigma,0}(\mathbb G)$ is 
$A(\mathbb G')$-$A(\mathbb G)$-biGalois.
The Hopf algebras $A(\mathbb G)$ and $A(\mathbb G')$
are isomorphic if and only if
there exists $u \in {\rm Aut}_g(G)$
such that $(u, \overline{\sigma}) \in \Gamma(\mathbb G)$.
\end{prop}

\noindent
\textbf{Proof}. It is easy to see that $\mathbb G'$ is a type IV 
group datum and it is straightforward to check that there exists
an algebra morphism $\beta : A_{\sigma, 0}(\mathbb G)
\longrightarrow A(\mathbb G') \otimes A_{\sigma, 0}(\mathbb G)$
such that $\beta(X) = 1 \otimes X + x \otimes T_g$
and $\beta(T_h) = h \otimes T_h$, $\forall h \in G$,   
that endows $A_{\sigma, 0}(\mathbb G)$ with an 
$A(\mathbb G')$-$A(\mathbb G)$ bicomodule algebra structure.
Very similarly to Proposition 3.3, we see that
$A_{\sigma, 0}(\mathbb G)$ is $A(\mathbb G')$-$A(\mathbb G)$-biGalois.
The Hopf algebras $A(\mathbb G')$ and $A(\mathbb G)$ are isomorphic
if and only if the corresponding group data are isomorphic, i.e.
if and only if there exists $u \in {\rm Aut}_g(G)$ such that
$\chi \circ u = \chi '$. This last condition exactly means that
$(u, \overline{\sigma}) \in \Gamma(\mathbb G)$. $\square$

\medskip

It is possible that the Hopf algebras $A(\mathbb G)$ and $A(\mathbb G')$
are not isomorphic. Indeed we have the following example.

\begin{ex}
{\rm 
Let $d >1$ and let $G = \langle x,h \ | \ x^{d^4} = 1 = h^{d^2} \ , \ 
xh = hx \rangle$ ($G \cong C_{d^4} \times C_{d^2}$).
Let $g = x^d$ and let $q\in \dot k$ be a primitive $d^2$th root of
unity. Let $\chi : G \longrightarrow \dot k$ be the character
defined by $\chi(x) = q$ and $\chi(h)=1$.
We have $d = o(\chi(g))$, $d^2 = o(\chi)$ and $d^3 = o(g)$.
For any $\sigma \in Z^2(G, \dot k)$, we have
$\sigma (g^d,x) =\sigma(x,g^d)$ while $\chi^d(x) = \chi(g) \not = 1$.
Hence $\mathbb G = (G,g, \chi)$ is a type IV group datum.

Now consider $\sigma \in Z^2(G,\dot k)$ defined by
$\sigma(x^\alpha h^\beta,x^{\alpha '}h^{\beta '}) = q^{\alpha \beta '}$
($\sigma$ is in fact a bicharacter).
We have $\sigma (g,h)^{-1} \sigma(h,g) = q^{-1}$.
Let $u : G \longrightarrow G$ be a group morphism. We 
have $u(h) = x^\alpha h^\beta$ and since $u(h^{d^2}) = 1$, we 
have $d^2 | \alpha$ and we conclude that $\chi\circ u(h) = 1$.
Hence there does not exist $u \in {\rm Aut}_g(G)$
such that $(u, \overline{\sigma}) \in \Gamma(\mathbb G)$, and 
we conclude that the Hopf algebras
$A(\mathbb G)$ and $A(\mathbb G')$ are not isomorphic.
}
\end{ex}

The following result concludes the proof of Theorem 3.2
in the type IV case. Since its proof is very similar to the 
one of Proposition 3.8, it is left to the reader.

\begin{prop}
Let $\mathbb G = (G,g, \chi)$ be a type IV group datum.
Then we have a group isomorphism:
\begin{align*}
\Psi : \Gamma(\mathbb G) 
\longrightarrow & {\rm BiGal}(A(\mathbb G)) \\
(u, \overline{\sigma}) \longmapsto & [A_{\sigma,0}^u].
\end{align*}
\end{prop}

We now study the case of type III group data. 
This is certainly the richest case.  

\begin{prop} 
Let $\mathbb G = (G,g, \chi)$ be a type III group datum.

\noindent
i) Let $\sigma \in Z^2(G, \dot k)$.
Consider the group morphism $\chi ' : G \longrightarrow \dot k$
defined by $\chi'(h) = \sigma(g,h)^{-1} \sigma(h,g) \chi(h)$
and the group datum $\mathbb G' = (G,g , \chi ')$.
The group datum $\mathbb G '$ is a type III or V group datum.
The right $A(\mathbb G)$-Galois object 
$A_{\sigma,0}(\mathbb G)$ is 
$A(\mathbb G')$-$A(\mathbb G)$-biGalois, and the Hopf 
algebras $A(\mathbb G)$ and $A(\mathbb G')$
are isomorphic if and only if
there exists $u \in {\rm Aut}_g(G)$
such that $(u, \overline{\sigma}) \in \Gamma(\mathbb G)$.

\noindent
ii) Let $\sigma \in Z^2(G, \dot k)$ and let $a \in \dot k$ such that
$A_{\sigma, a}(\mathbb G)$ is $A(\mathbb G)$-Galois.
Consider the group morphism $\chi ' : G \longrightarrow \dot k$
defined by $\chi'(h) = \sigma(g,h)^{-1} \sigma(h,g) \chi(h)$
and let $\mu = - a \sigma(g,g)^{-1} \ldots \sigma(g,g^{d-1})^{-1}$.
Then $\mathbb G' = (G,g,\chi', \mu)$ is a type VI group datum,
and $A_{\sigma, a}(\mathbb G)$ is $A(\mathbb G')$-$A(\mathbb G)$-biGalois.

\end{prop}

\noindent
\textbf{Proof}. i) The proof is similar to the one of Proposition 3.10.

\noindent
ii) Since $A_{\sigma, a}(\mathbb G)$ is $A(\mathbb G)$-Galois,
we have $\sigma(g^d,h)= \sigma(h,g^d)$, $\forall h \in G$. 
Hence $\chi'^d = 1$ and $\mathbb G'$ is a type VI group datum. 
Then one checks that there exists
an algebra morphism $\beta : A_{\sigma, a}(\mathbb G)
\longrightarrow A(\mathbb G') \otimes A_{\sigma, a}(\mathbb G)$
such that $\beta(X) = 1 \otimes X + x \otimes T_g$
and $\beta(T_h) = h \otimes T_h$, $\forall h \in G$,   
that endows $A_{\sigma, a}(\mathbb G)$ with an 
$A(\mathbb G')$-$A(\mathbb G)$ bicomodule algebra structure.
Very similarly to Proposition 3.3, we see that
$A_{\sigma, a}(\mathbb G)$ is $A(\mathbb G')$-$A(\mathbb G)$-biGalois.
$\square$

\medskip

The following two examples show that
the group datum $\mathbb G'$ of proposition 3.13.i might
be of type III not isomorphic with $\mathbb G$, or of type V.

\begin{ex}
{\rm 
Let $d >1$ and let $G = \langle g,h \ | \ g^{d^2} = 1 = h^{d} \ , \ 
gh = hg \rangle$ ($G \cong C_{d^2} \times C_{d}$).
Let $q\in \dot k$ be a primitive $d$th root of
unity. Let $\chi : G \longrightarrow \dot k$ be the character
defined by $\chi(g) = q$ and $\chi(h)=1$: we have $g^d \not = 1$
and $\chi^d = 1$ and hence 
$\mathbb G = (G,g, \chi)$ is a type III group datum.

Now consider $\sigma \in Z^2(G,\dot k)$ defined by
$\sigma(g^\alpha h^\beta,g^{\alpha '}h^{\beta '}) = q^{\alpha \beta '}$
($\sigma$ is in fact a bicharacter).
We have $\sigma (g,h)^{-1} \sigma(h,g) = q^{-1}$ and 
hence $\chi'^d =1$: $\mathbb G' = (G,g,\chi')$ is type III group datum.
Let $u : G \longrightarrow G$ be a group morphism. We 
have $u(h) = g^\alpha h^\beta$ and since $u(h^{d}) = 1$, we 
have $d | \alpha$ and we conclude that $\chi\circ u(h) = 1
\not = q^{-1}$.
Hence there does not exist $u \in {\rm Aut}_g(G)$
such that $(u, \overline{\sigma}) \in \Gamma(\mathbb G)$, and 
we conclude that the Hopf algebras
$A(\mathbb G)$ and $A(\mathbb G')$ are not isomorphic.
}
\end{ex}

\begin{ex}
{\rm 
Let $d >1$ and let $G = \langle g,h \ | \ g^{d^2} = 1 = h^{d^2} \ , \ 
gh = hg \rangle$ ($G \cong C_{d^2} \times C_{d^2}$).
Let $q\in \dot k$ be a primitive $d^2$th root of
unity. Let $\chi : G \longrightarrow \dot k$ be the character
defined by $\chi(g) = q^d$ and $\chi(h)=1$: we have $g^d \not = 1$
and $\chi^d = 1$ and hence 
$\mathbb G = (G,g, \chi)$ is a type III group datum.

Now consider $\sigma \in Z^2(G,\dot k)$ defined by
$\sigma(g^\alpha h^\beta,
g^{\alpha '}h^{\beta '}) = q^{\alpha \beta '}$
($\sigma$ is a bicharacter).
We have $\sigma (g,h)^{-1} \sigma(h,g) = q^{-1}$ and 
hence $\chi'^d (h) = \sigma(g^d,h)^{-1} \sigma(h,g^d) \chi^d(h)
= \sigma(g^d,h)^{-1} \sigma(h,g^d) = q^{-d}$. 
Hence $\mathbb G' = (G,g,\chi')$ is type V group datum.
}
\end{ex}

Using Proposition 3.13, the proof of the description
of the biGalois group in the type III case is the same has 
in the previous cases, and hence is left once again to the reader.

\begin{prop}
Let $\mathbb G = (G,g, \chi)$ be a type III group datum.
Then we have a group isomorphism:
\begin{align*}
\Psi : \Gamma(\mathbb G) 
\longrightarrow & {\rm BiGal}(A(\mathbb G)) \\
(u, \overline{\sigma}) \longmapsto & [A_{\sigma,0}^u].
\end{align*}
\end{prop}

\begin{coro}
Let $\mathbb G = (G,g, \chi)$ be a type V group datum, and let
$\sigma \in Z^2(G,\dot k)$ be such that
$\chi^d(h) = \sigma(g^d,h)\sigma(h,g^d)^{-1}$.
Consider the group morphism $\chi_{\sigma} : G \rightarrow \dot k$
defined by $\chi_{\sigma}(h) = \sigma(g,h)^{-1}\sigma(h,g)\chi(h)$.
Then $\mathbb G_{\sigma} = (G,g,\chi_{\sigma})$ is a type III
group datum, there exists an $A(\mathbb G)$-$A(\mathbb G_{\sigma})$-biGalois
extension and group isomorphisms
$${\rm BiGal}(A(\mathbb G))  \cong {\rm BiGal}(A(\mathbb G_{\sigma}))
\cong \Gamma(\mathbb G_{\sigma}) \cong \Gamma(\mathbb G).$$
\end{coro}

\noindent
\textbf{Proof}. It is immediate that $\mathbb G_{\sigma}$ is a type
III group datum. Consider now $\sigma^{-1} \in Z^2(G, \dot k)$.
It is clear that the associated group datum $(\mathbb G_{\sigma})'$ of
Proposition 3.13.i is $\mathbb G$, and hence that 
$A_{\sigma^{-1},0}(\mathbb G)$ is 
$A(\mathbb G)$-$A(\mathbb G_{\sigma})$-biGalois, which
gives the first isomorphism. The second isomorphism
is given by Proposition 3.16. For the last one,
it is not difficult to check that we have a map
\begin{align*}
\Gamma(\mathbb G_{\sigma}) \longrightarrow & \Gamma(\mathbb G) \\
(u, \overline{\tau}) \longmapsto & 
(u, \overline{(\sigma\circ u \times u)^{-1} \sigma \tau})
\end{align*}
which is a group isomorphism.
$\square$

\medskip

This last corollary finally completes the proofs
of Theorems 2.1 and 3.2.
              
\begin{coro}
Let $\mathbb G = (G,g, \chi, \mu)$ be a type VI group datum.
Consider the type III group datum 
$\mathbb G_{red} = (G,g,\chi)$. Then there exists
an $A(\mathbb G)$-$A(\mathbb G_{red})$ biGalois extension, and hence
a bijection
$${\rm Gal}(A(\mathbb G)) \cong H^2(G, \dot k) \amalg
H^2_{g^d,g^d}(G, \dot k)$$
and a group isomorphism
$${\rm BiGal}(A(\mathbb G)) \cong \Gamma(\mathbb G).$$
\end{coro}

\noindent
\textbf{Proof}. It is clear from the second
part of Proposition 3.13 that
$A_{1, -\mu}(\mathbb G_{red})$ is 
$A(\mathbb G)$-$A(\mathbb G_{red})$-biGalois. Hence the type III
case of Theorems 2.1 and 3.2 conclude the proof. $\square$

\section{Examples}

This section is devoted to the study of a few
concrete examples. We first consider cyclic group data, 
and we get in particular Masuoka's description \cite{[Ma]}
of the Galois objects for the Taft algebras and 
Schauenburg's description \cite{[Sc1]} of their biGalois groups.
Then we consider decomposable group data and
get a more explicit description of corresponding
biGalois groups: in particular we recover Schauenburg's 
description of the biGalois groups of the generalized
Taft algebras, in a slightly different form.

\subsection{Cyclic group data}

A group datum $\mathbb G = (G,g,\chi,\mu)$ is said to be \textbf{cyclic} 
if the underlying group $G$ is cyclic. In this subsection
we completely describe the Galois objects and biGalois groups
for Hopf algebras associated with cyclic group data.

Before stating the main result, we need to introduce some terminology. 
Let $N > 1$ be an integer and let $n_1, \ldots , n_r$ be some divisors
of $N$. We put
$$U(\mathbb Z / N \mathbb Z)[n_1, \ldots , n_r] =
\{ \overline{\beta} \in U(\mathbb Z/ N \mathbb Z)
\ (\beta \in \mathbb Z) \ | \ \beta \equiv 1 \ ({\rm mod} \ n_i), \
i=1, \ldots , r \}.$$  
This is clearly a subgroup of $U(\mathbb Z / N \mathbb Z) \cong
{\rm Aut}(C_N)$. For $N=1$, we adopt the convention
$U(\mathbb Z / N \mathbb Z) = \{1\}$.

\begin{theo}
Let  $\mathbb G = (G,g,\chi,\mu)$ be a cyclic group datum
with $d= o(\chi(g))$, $n= o(g)$, $N = |G|$ and $m = o(\chi)$.

\noindent
$\bullet$
If $\mathbb G$ is a type I group datum, i.e. if $d=n=m$, then
$${\rm Gal}(A(\mathbb G)) \cong \dot k/ \dot k^N \times k
\quad {\rm and} \quad
{\rm BiGal}(A(\mathbb G)) \cong
({\rm Aut}(C_{\frac{N}{n}}) 
\ltimes  (\dot k \times \dot k/ \dot k^{\frac{N}{n}})) \ltimes k.$$
$\bullet$ 
If $\mathbb G$ is a type II group datum, i.e. if $d=n<m$, then
$${\rm Gal}(A(\mathbb G)) \cong \dot k/ \dot k^N 
\quad {\rm and} \quad
{\rm BiGal}(A(\mathbb G)) \cong
U(\mathbb Z / N \mathbb Z)[m] 
\ltimes  (\dot k \times \dot k/ \dot k^{\frac{N}{n}}).$$
$\bullet$
If $\mathbb G$ is a type III or VI group datum, i.e. if $d=m<n$, then
$${\rm Gal}(A(\mathbb G)) \cong \dot k/ \dot k^N \amalg
(\dot k \times  \dot k/ \dot k^{\frac{Nd}{n}})
\quad {\rm and} \quad
{\rm BiGal}(A(\mathbb G)) \cong
U(\mathbb Z / N \mathbb Z)[n]
\ltimes  (\dot k \times \dot k/ \dot k^{\frac{N}{n}}).$$
$\bullet$
If $\mathbb G$ is a type IV group datum, i.e. if $d<n$ and $d<m$, then
$${\rm Gal}(A(\mathbb G)) \cong \dot k/ \dot k^N 
\quad {\rm and} \quad
{\rm BiGal}(A(\mathbb G)) \cong
U(\mathbb Z / N \mathbb Z)[n,m]
\ltimes  (\dot k \times \dot k/ \dot k^{\frac{N}{n}}).$$
\end{theo} 

The rest of this subsection is essentially devoted to the
proof of Theorem 4.1. The first remark is that since any 2-cocycle 
on a cyclic group is symmetric, there does not exist any type V
cyclic group datum, and we have the following result.

\begin{lemm}
Let $\mathbb G = (G,g,\chi,\mu)$ be a cyclic group
datum. Then $$\Gamma(\mathbb G) \cong 
{\rm Aut}_{g,\chi}(G) \ltimes H^2_{1,g}(G, \dot k).$$
\end{lemm}

The next step is to compute the modified cohomology groups
of the first section in the case of a cyclic group.
This is done in the following lemma.

\begin{lemm} Let $C_N$ be the cyclic group of order $N$ and let
$g \in C_N$ be an element of order $n>1$. Then we have a group
isomorphism:
$$H^2_{1,g}(C_N,\dot k) = H_{g,g}^2(C_N,\dot k) 
\cong \dot k \times ( \dot k / \dot k^{\frac{N}{n}}).$$
\end{lemm}

\noindent
\textbf{Proof}. 
We construct a split exact sequence
$$
\begin{CD}
1 @>>> \dot k / \dot k^{\frac{N}{n}} @>u>> H^2_{1,g}(G, \dot k)
 @>\epsilon>> \dot k @>>> 1
\end{CD}
$$
where the group morphism
$\epsilon : H^2_{1,g}(G,\dot k) \longrightarrow \dot k$ was constructed
in subsection 1.3.
Let us choose a generator $y$ of $C_N$
such that  $g=y^{\frac{N}{n}}$. 
For $a \in \dot k$, define $f_a \in Z^2(G,\dot k)$ 
in the following way \cite{[Kar]}. For
$0 \leq i,j \leq  N-1$, we put
$f_a(y^i,y^j) = 1$ if $i+j \leq N-1$ and
$f_a(y^i,y^j)=a$ if $i+j \geq N$.
We have $\epsilon(\overline{f_a}) = a$ and 
hence the group morphism $f : \dot k \rightarrow H^2_{1,g}(G,\dot k)$,
$a \longmapsto \overline{f_a}$, satisfies $\epsilon \circ f = 
{\rm id}_{\dot k}$.

Let us fix a family of integers
$(r_i)_{0 \leq i \leq N-1}$ with $r_0=1$ and $r_{\frac{N}{n}} = 1$, and 
for $t \in \dot k$, let us define $\beta_t : G \longrightarrow \dot k$ by
$\beta_t(y^i) = t^{r_i}$, for $0 \leq i \leq N-1$. 
We have $\beta_t(g) =t$.  Now define a group morphism
$$
u_0 :  \dot k  \longrightarrow H_{1,g}^2(G,\dot k), \quad
 t \longmapsto \overline{f_{t^{-n}} \partial (\beta_t)}.$$
It is clear that $u_0(\dot k) \subset {\rm Ker} (\epsilon)$.
Now let $\sigma \in Z^2(G, \dot k)$ be such that
$\overline{\sigma} \in {\rm Ker}(\epsilon)$.
By the description of $H^2(G,\dot k)$ in \cite{[Kar]},
Theorem 2.3.1, there exists
$a \in \dot k$ and $\mu : G \rightarrow \dot k$ (with $\mu(1)=1$)
such that $\sigma = f_a \partial (\mu)$.  
Then $a = \mu(g)^{-n}$. Since $\beta_{\mu(g)}(g) = \mu(g)$,
we have $\overline{\partial(\mu)} = \overline{\partial(\beta_{\mu(g)}})$,
and $\overline{\sigma} = \overline{f_{\mu(g)^{-n}}\partial(\beta_{\mu(g)})}
\in u(\dot k)$. We conclude that 
Ker$(\epsilon) = u_0(\dot k)$. 

Now let $t \in \dot k$ be such that
$u_0(t) = \overline{1}$. Then there exists $\mu : G \rightarrow \dot k$
with $\mu(g)=1$ such that 
$f_{t^{-n}} \partial(\beta_t)= \partial(\mu)$, and hence 
$f_{t^{-n}}$ is a coboundary, which means that $t^{-n} \in \dot k^N$
(\cite{[Kar]}), i.e. that $t \in \dot k^{\frac{N}{n}}$
(recall that $k$ contains all primitive roots of unity). 
Conversely let $t \in \dot k^{\frac{N}{n}}$ and let $s \in\dot k$
be such that $t^{-n} = s^N$. Define $\gamma_s : G \longrightarrow \dot k$ by
$\gamma_s(y^i) =s^i$ for $0 \leq i \leq N-1$.
Then $f_{t^{-n}} = \partial(\gamma_s)$ and 
$f_{t^{-n}} \partial(\beta_t) = \partial(\gamma_s \beta_t)$,
with $\gamma_s \beta_t(g) = s^{\frac{N}{n}}t =1$, and thus
$t \in {\rm Ker}(u_0)$. Therefore
$u_0$ induces an injective morphism
$u : \dot k/ \dot k^{\frac{N}{n}} \longrightarrow H^2_{1,g}(G,\dot k)$
and we have the announced split exact sequence. $\square$   

\medskip

Let us now describe the explicit models for cyclic group
data (we do not treat the type VI case thanks to Corollary 
3.18). A \textbf{cyclic datum} is a $5$-uplet
$(d,n,N,\alpha,q)$ where $d,n,N>1$ are integers, $\alpha \in \mathbb N^*$
and $q \in \dot k$ is a root of unity, satisfying:
$$d|n|N, \ \alpha | \frac{N}{n}, \ {\rm GCD}(\alpha,d)=1, \ o(q) = 
\frac{Nd}{\alpha n}.$$ 
To any cyclic datum $(d,n,N,\alpha,q)$
we associate a group datum
$$\mathbb C[d,n,N,\alpha,q] :=
(C_N=\langle z \rangle, g = z^{\frac{N}{n}}, \chi_q)$$
where $\chi_q$ is the character defined by $\chi_q(z) = q$. 
The exact relations between cyclic data and cyclic group data
are given by the following lemma.

\begin{lemm}
i) Let $(d,n,N,\alpha,q)$ be a cyclic datum.
Then $\mathbb C[d,n,N,\alpha,q]$ is a cyclic group
datum with $d= o(\chi_q(g))$, $n=o(g)$, $N=|G|$
and $m = o(\chi_q) = \frac{Nd}{\alpha n}$.

\noindent 
ii) Let $\mathbb G=(G,g,\chi)$ be a cyclic group datum. 
Then there exists a cyclic datum $(d,n,N,\alpha,q)$
such that $\mathbb G \cong \mathbb C[d,n,N,\alpha,q]$.

\noindent
iii) Let $(d,n,N,\alpha,q)$ be a cyclic datum, with 
the cyclic group datum $\mathbb C[d,n,N,\alpha,q]$.

\noindent
$\bullet$ $\mathbb C[d,n,N,\alpha,q]$ is a type I group datum
if and only if
$$d= N \quad {\rm or} \quad [d=n<N , \ {\rm GCD}(\frac{N}{n},n) =1 \
{\rm and} \ \alpha = \frac{N}{n}].$$
$\bullet$ $\mathbb C[d,n,N,\alpha,q]$ is a type II group datum
if and only if
$$d= n < N \ {\rm and} \ \alpha < \frac{N}{n}.$$
$\bullet$ $\mathbb C[d,n,N,\alpha,q]$ is a type III group datum
if and only if
$$d<n= N \quad {\rm or} \quad [d<n<N, \ {\rm GCD}(\frac{N}{n},d) =1 \
{\rm and} \ \alpha = \frac{N}{n}].$$
$\bullet$ $\mathbb C[d,n,N,\alpha,q]$ is a type IV group datum
if and only if
$$d< n<N \ {\rm and} \ \alpha < \frac{N}{n}.$$
\end{lemm}

\noindent
\textbf{Proof}. i) We have 
$$o(\chi(g)) = o(q^{\frac{N}{n}})= 
\frac{o(q)}{{\rm GCD}(o(q),\frac{N}{n})}=
\frac{\frac{Nd}{\alpha n}}{{\rm GCD}({\frac{N}{\alpha n}\alpha,
\frac{N}{\alpha n} d})}=d$$
since GCD$(\alpha,d)=1$. The other assertions are immediate.

\noindent
ii) Let $\mathbb G = (G,g, \chi)$ be a cyclic group datum
with $|G|=N$ and $o(g)=n$ 
Let us choose a generator $z \in G$ such that $g = z^{\frac{N}{n}}$.
Let $q = \chi(z)$. We have 
$$d=o(\chi(g))=  o(q^{\frac{N}{n}}) = 
\frac{o(q)}{{\rm GCD}(o(q),\frac{N}{n})}.$$ 
From this we see that $o(q)|(\frac{N}{n}d)$.
Let $\alpha \in \mathbb N^*$ be such that
$\alpha o(q)= \frac{N}{n} d$. Then
$\frac{o(q)}{d} \alpha= \frac{N}{n}$ and hence
$\alpha | \frac{N}{n}$. Finally
$$\frac{N}{n \alpha} = {\rm GCD}(o(q),\frac{N}{n})
={\rm GCD}(\frac{N}{n\alpha}d,\frac{N}{n\alpha} \alpha)$$
and we conclude that ${\rm GCD}(\alpha,d)=1$. 
Thus $(d,n,N, \alpha, q)$ is a cyclic datum and we
have $\mathbb G \cong \mathbb C[d,n,N,\alpha,q]$.
The final assertions are immediate. $\square$

\medskip

\noindent
\textbf{Proof of Theorem 4.1}. Let $\mathbb G = (G,g,\chi)$ be a
cyclic group datum. By Lemma 4.4 we can assume that
$$\mathbb G = \mathbb C[d,n,N,\alpha,q] =
(C_N=\langle z \rangle, g = z^{\frac{N}{n}}, \chi_q)$$
for a cyclic datum $(d,n,N,\alpha, q)$. 
We begin with some generalities.
We always have $H^2_{1,g}(G,\dot k) \cong \dot k \times 
\dot k / \dot k^{\frac{N}{n}}$ by Lemma 4.3. 
Let $u \in {\rm Aut}_{g,\chi}(G)$. Then there exists
$\beta \in \{0, \ldots ,N-1\}$ such that $u(z) =z^\beta$ with 
$${\rm GCD}(\beta, N)=1, \quad n | (\beta -1), \quad
\frac{Nd}{\alpha n} | (\beta-1).$$
Hence we have
$${\rm Aut}_{g,\chi}(G) \cong 
U(\mathbb Z/N\mathbb Z)[n, \frac{Nd}{\alpha n}]
\quad {\rm and} \quad
\Gamma(\mathbb G) \cong U(\mathbb Z/N\mathbb Z)[n, \frac{Nd}{\alpha n}]
\ltimes  (\dot k \times \dot k / \dot k^{\frac{N}{n}})$$
by Lemma 4.2. The proof of Theorem 4.1 now follows easily
from Theorems 2.1, 3.2, Lemma 4.4 and the case by case computation
of $U(\mathbb Z/N\mathbb Z)[n, \frac{Nd}{\alpha n}]$. $\square$

\medskip

Let us specialize Theorem 4.1 to a few famous examples
of Hopf algebras. First consider the cyclic datum
$(N,N,N,1,q)$. The cyclic group datum
$\mathbb C[N,N,N,1,q]$ is of type I, and the corresponding
Hopf algebra $A(\mathbb C[N,N,N,1,q])$ is precisely the Taft
algebra $H_{N,q}$. Hence we 
get Masuaoka's \cite{[Ma]} and Schauenburg's \cite{[Sc2]}
respective results 
$${\rm Gal}(H_{N,q}) \cong \dot k / \dot k^{\frac{N}{n}}
\quad {\rm and} \quad
{\rm BiGal}(H_{N,q}) \cong \dot k \ltimes k.$$
Now consider the cyclic datum
$(d,N,N,1,q)$ with $d<N$. The cyclic group datum 
$\mathbb C[d,N,N,1,q]$ is of type III, and hence
for $\mu \in k$ we can consider the group datum
$\mathbb C[d,N,N,1,q][\mu]
= (G = C_N = \langle g \rangle, g , \chi_q,\mu)$. 
The corresponding Hopf algebra
$A(\mathbb C[d,N,N,1,q][\mu])$ is,  with the notation
of \cite{[Ra]}, the simple-pointed Hopf algebra
$A_{(q,\mu,d,N)}$ of \cite{[Ra]}. Hence we get
$${\rm Gal}(A_{(q,\mu,d,N)}) \cong \dot k / \dot k^N
\amalg (\dot k \times k / \dot k^d)
\quad {\rm and} \quad
{\rm BiGal}(A_{(q,\mu,d,N)}) \cong \dot k.$$  

\subsection{Decomposable group data}

The first goal of this section is to provide a more concrete
description for the biGalois group of the Hopf algebra
associated to a type I group datum. We prove the following result.

\begin{theo}
Let $\mathbb G = (G,g,\chi)$ be a type I group datum
with $K = {\rm Ker}(\chi)$ and $d=o(\chi(g)) =o(g) =o(\chi)$. 
Then we have a bijection
$${\rm Gal}(A(\mathbb G)) \cong
\dot k/ \dot k^d \times H^2(K,\dot k) \times {\rm Hom}(K,\mu_d) \times k,$$
and a group isomorphism
$${\rm BiGal}(A(\mathbb G)) \cong {\rm Aut}(K) \ltimes
((\dot k \ltimes k) \times H^2(K,\dot k) \times
{\rm Hom}(K,\mu_d))$$
where the right actions of ${\rm Aut}(K)$ on $H^2(K, \dot k)$
and ${\rm Hom}(K, \mu_d)$ are the natural ones, and the
action on $\dot k \ltimes k$ is the trivial one.
\end{theo}

In fact the comparison of 
our description of the biGalois group of a generalized Taft algebra
with Schauenburg's one \cite{[Sc2]} will lead us to consider
more generally decomposable group data. 
We say that a  group datum $\mathbb G =(G,g,\chi,\mu)$
is \textbf{decomposable} if $G$ is the direct product of $\langle g \rangle$
with a subgroup $K$. It is easy to see that $\mathbb G$ is a type
I group datum if and only if $\mu=0$ and $\mathbb G$ is decomposable
with $G$ being the direct product of $\langle g \rangle$
and of $K = {\rm ker}(\chi)$.

The first step when studying decomposable group data is to
examine the behaviour of the modified cohomology group
under direct product.
We have the following slight generalization of a theorem of 
Yamazaki, see \cite{[Kar]}, Theorem 2.3.13.
For groups $G_1$ and $G_2$, the set of pairings (i.e. bimorphisms)
$G_1 \times G_2 \longrightarrow \dot k$ is denoted
by $\mathcal P(G_1 \times G_2, \dot k)$. It has a natural group
structure.

\begin{lemm}
Let $G_1,G_2$ be some groups and let $g,h \in Z(G_1)$.
Then we have a group isomorphism
$$H^2_{g,h}(G_1 \times G_2, \dot k) \cong 
H^2_{g,h}(G_1, \dot k) \times H^2(G_2, \dot k)
\times \mathcal P(G_1/\langle g \rangle \times G_2,\dot k).$$
\end{lemm}

\noindent
\textbf{Proof}. We freely use the techniques and results of
the proof of Theorem 2.3.13 in \cite{[Kar]}.
Let us first consider the group morphism
\begin{align*}
\Psi_0 : 
Z^2_g(G_1  \times G_2, \dot k) \longrightarrow &
Z^2_g(G_1  ,\dot k) \times Z^2(G_2,\dot k) \times
\mathcal P(G_1/\langle g \rangle \times G_2,\dot k) \\ 
 \sigma \longmapsto & (\sigma_1, \sigma_2, B_{\sigma})
\end{align*}
where $\sigma_i$, $i=1,2$, is the restriction of $\sigma$ to $G_i$
and $B_{\sigma}$ is the pairing defined by
$B_{\sigma}(\overline{g_1},g_2) = \sigma(g_1,g_2) \sigma(g_2,g_1)^{-1}$,
$\forall (g_1,g_2) \in G_1 \times G_2$.
Clearly $\Psi_0$ induces a group morphism
\begin{align*}
\Psi : 
H^2_{g,h}(G_1  \times G_2, \dot k) \longrightarrow &
H^2_{g,h}(G_1  ,\dot k) \times H^2(G_2,\dot k) \times
\mathcal P(G_1 / \langle g  \rangle \times G_2,\dot k) \\ 
\overline{ \sigma} \longmapsto & 
(\overline{\sigma_1}, \overline{\sigma_2}, B_{\sigma}).
\end{align*} 
For $(\sigma_1, \sigma_2, B) \in  Z^2_g(G_1  ,\dot k) 
\times Z^2(G_2,\dot k) \times \mathcal P(G_1/\langle g \rangle 
\times G_2,\dot k)$, define
an element $\sigma \in Z^2_g(G_1 \times G_2,\dot k)$ by
$$\sigma(g_1g_2,g'_1g'_2) = 
\sigma_1(g_1,g'_1)\sigma_2(g_2,g'_2) B(\overline{g_1},g'_2), \ 
\forall g_1,g'_1 \in G_1, \ g_2,g'_2 \in G_2.$$
Clearly $\Psi(\overline{\sigma}) =
(\overline{\sigma_1}, \overline{\sigma_2}, B_{\sigma})$, and
hence $\Psi$ is surjective.
Now let $\sigma \in Z^2_g(G_1 \times G_2,\dot k)$ be such that
$\overline{\sigma} \in {\rm Ker}(\Psi)$. 
This means that there exists $\mu_1 : G_1 \rightarrow \dot k$
with $\mu_1(h) =1 = \mu_1(1)$, $\mu_2 : G_2 \rightarrow \dot k$
with $\mu_2(1)=1$, and that $B_\sigma= 1$.
Now define a map 
$$\mu : G_1 \times G_2 \rightarrow \dot k \ {\rm by} \
\mu(g_1g_2) = \mu_1(g_1)\mu_2(g_2), \ \forall (g_1,g_2) \in G_1 \times G_2.$$
We have $\mu(h)=1$ and hence 
$\overline{\sigma} = \overline{\partial(\mu) \sigma}$.
In this way we can assume without loss of generality that
$\sigma(g_1,g'_1) = 1 = \sigma(g_2,g'_2)$, 
$\forall g_1,g'_1 \in G_1, \ g_2,g'_2 \in G_2$.
Now consider the map $s : G_1 \times G_2 \rightarrow \dot k$ defined
by $s(g_1,g_2) = \sigma(g_1,g_2)^{-1}$. By the computation
in the proof of Theorem 2.3.13 in \cite{[Kar]} (p. 61), we have
$\sigma = \partial (s)$, and since $s(h)=1$, we conclude
that $\overline{\sigma} = \overline{1}$ and that 
$\Psi$ is injective. $\square$ 

\medskip

Using this lemma, or simply Yamazaki's result 
in the type I, II and IV cases, one gets a more precise description
for the set ${\rm Gal}(A(\mathbb G))$ of a decomposable group
datum. The type I case is contained in the statement
of Theorem 4.5, and we have:

\begin{coro}
Let $\mathbb G = (G,g,\chi)$ be a decomposable group datum
with $G = \langle g \rangle \times K$, and $d=o(\chi(g))$, $n = o(g)$.
Assume that $\mathbb G$ is not a type I group datum.

\noindent
$\bullet$ If $\mathbb G$ is of type II or IV, then
$${\rm Gal}(A(\mathbb G)) \cong
\dot k/ \dot k^n \times H^2(K,\dot k) \times {\rm Hom}(K,\mu_n)$$

\noindent
$\bullet$
If $\mathbb G$ is of type III, V or VI, then
$${\rm Gal}(A(\mathbb G)) \cong
(\dot k/ \dot k^n \times H^2(K,\dot k) \times {\rm Hom}(K,\mu_n))
\amalg (\dot k \times 
\dot k/ \dot k^d \times H^2(K,\dot k) \times {\rm Hom}(K,\mu_n)).$$
\end{coro}

We now study the biGalois group of a type I group datum.

\begin{lemm}
Let $\mathbb G = (G,g,\chi)$ be a type I group datum, let $K = {\rm Ker}
(\chi)$ and let $p_2 : G \longrightarrow K$ be the canonical
projection. The map
\begin{align*}
\Omega : \Gamma(\mathbb G) \longrightarrow & {\rm Aut}(K) \ltimes
(\dot k  \times H^2(K,\dot k) \times {\rm Hom}(K,\mu_d)) \\
(u, \overline{\sigma}) \longmapsto &
(p_2 \circ u_{|K}, \sigma(g,g) \ldots \sigma(g,g^{d-1}) ,
\overline{\sigma_{|K\times K}}, \chi  \circ u_{|K})
\end{align*}
is a group isomorphism.
\end{lemm}

\noindent
\textbf{Proof}. It is straightforward to check that
$\Omega$ is well defined (is fact on the whole group
Aut$_g(G) \ltimes H^2_{1,g}(G, \dot k)$) and is a group
morphism. Let $(u,\overline{\sigma}) \in {\rm Ker \Omega}$.
It is easy to see that, since $G = \langle g\rangle K$ and 
$K={\rm Ker}(\chi)$, that $u={\rm id}_G$.
Since $(u,\overline{\sigma}) \in  \Gamma(\mathbb G)$, we
have $\chi \circ u(h) = \sigma(g,h)^{-1} \sigma(h,g) \chi(h)$,
$\forall h \in G$, and hence $\sigma(g,h) = \sigma(h,g)$, $\forall h \in G$.
Thus the pairing $B_{\sigma}$ in the proof of Lemma 4.6 is trivial.
By Lemma 4.3 the identity $\sigma(g,g) \ldots \sigma(g,g^{d-1})=1$
precisely means that the restriction of $\sigma$ to $\langle g \rangle
\times \langle g \rangle$ is trivial in $H^2_{1,g}(G,\dot k)$, and we conclude
that $\overline{\sigma} =\overline{1}$ by Lemma 4.6:
$\Omega$ is injective. 

Finally let $(f, \lambda, \overline{\tau}, \psi) \in 
{\rm Aut}(K) \times
\dot k  \times H^2(K,\dot k) \times {\rm Hom}(K,\mu_d)$.
Let $s : \mu_d \rightarrow \langle g \rangle$ be a group morphism
such that $\chi \circ s = {\rm id}$. Define
$u \in {\rm Aut}_g(G)$ by $u(g) = g$
and $u(h) =s(\psi(h)) f(h)$ for $h \in K$. 
Let $f_\lambda \in Z^2(\langle g \rangle, \dot k)$ defined
in the proof of Lemma 4.3, and define $\sigma \in Z^2(G, \dot k)$
as in the proof of Lemma 4.6:
$$\sigma(g^\alpha h,g^\beta h') = 
f_\lambda(g^\alpha,g^\beta)\tau(h,h') \psi(h')^\alpha, \ 
\forall \alpha,\beta \in \mathbb Z, \forall \ h,h' \in K.$$
It is clear that $(u,\overline{\sigma}) \in \Gamma(\mathbb G)$
and that $\Omega(u,\overline{\sigma}) = (f, \lambda, \overline{\tau}, \psi)$.
Thus $\Omega$ is an isomorphism. $\square$

\medskip

Theorem 4.5 is now a straightforward consequence of Theorem 3.2
and Lemma 4.8. $\square$

\medskip 

Before coming to the generalized Taft algebras, let us examine
the case of a type I group datum $\mathbb G = (G,g, \chi)$
where $K= {\rm Ker} (\chi)$ is an abelian group, written as 
a direct product of cyclic groups $K = \prod_{i=1}^m C_{N_i}$.
We have a group isomorphism
$$H^2(K , \dot k) \cong \prod_{i=1}^m H^2(C_{N_i}, \dot k) \times
{\rm Hom}(\Lambda^2(K),\dot k),$$ 
where ${\rm Hom}(\Lambda^2(K),\dot k)$ is the group
of alternating bimorphism $B : K \times K \rightarrow \dot k$
(i.e. with $B(r,r)=1$ for $r \in K$). 
This isomorphism is Aut$(K)$-equivariant,
the action of Aut$(K)$ on Hom$(\Lambda^2(K),\dot k)$ being the natural
one while we need a few words to explain the action
on $\prod_{i=1}^m H^2(C_{N_i},\dot k)$. Let 
$\underline{\sigma} =(\overline{\sigma_i}) 
\in \prod_{i=1}^m H^2(C_{N_i},\dot k)$ and let $u \in {\rm Aut}(K)$
written in a matrix form $u=(u_{ij})$ with $u_{ij} \in {\rm Hom}(G_j,G_i)$.
Then 
$$\underline{\sigma} \leftarrow u = 
(\prod_j \overline{\sigma_j \circ (u_{ji} \times u_{ji})})_i.$$
Thus by Theorem 4.5 we get a bijection
$$
{\rm Gal}(A(\mathbb G)) \cong
\dot k / \dot k^d  \times (\prod_{i=1}^m \dot k / \dot k^{N_i}) \times
{\rm Hom}(\Lambda^2(K),\dot k) \times {\rm Hom}(K,\mu_d) \times k,$$
and a group isomorphism
$${\rm BiGal}(A(\mathbb G)) \cong {\rm Aut}(K) \ltimes
((\dot k \ltimes k) \times 
 (\prod_{i=1}^m \dot k / \dot k^{N_i}) \times
{\rm Hom}(\Lambda^2(K),\dot k) \times
{\rm Hom}(K,\mu_d)).$$
The generalized Taft algebra $H_{N,m+1,q}$ with $(m+1)$ grouplike
($m \in \mathbb N$) may be defined as the Hopf algebra
associated with the group datum
$\mathbb G = (C_N^{m+1},g,\chi_q)$ where $g$ is a generator
of the first cyclic copy of $C_N^{m+1}$, $q$ is a primitive
$N$th root of unity and $\chi_q$ is the character defined
by $\chi_q(g) = q$ and Ker$(\chi) = C_N^m$. We get, after
standard identifications, the result of Schauenburg \cite{[Sc2]}
$$
{\rm Gal}(H_{N,m+1,q}) \cong
\dot k / \dot k^N  \times (\dot k / \dot k^{N})^m \times
{\rm Skew}_m(\mathbb Z / N \mathbb Z) \times (\mathbb Z / N \mathbb Z)^m
\times k,$$
(where Skew$_m(\mathbb Z / N \mathbb Z)$ is the set of 
$\mathbb Z / N \mathbb Z$-valued skew-symmetric matrices
$r = (r_{ij})$, i.e. with $r_{ij} = -r_{ji}$ and $r_{ii} = 0$)
 and a group isomorphism
$${\rm BiGal}(H_{N,m+1,q}) \cong {\rm GL}_m(\mathbb Z / N \mathbb Z) \ltimes
((\dot k \ltimes k) \times 
 (\dot k / \dot k^{N})^m \times
{\rm Skew}_m(\mathbb Z / N \mathbb Z) \times
(\mathbb Z / N\mathbb Z)^m),$$
The right actions of GL$_m(\mathbb Z / N \mathbb Z)$ on each copy
of the right side being the following ones.
The action on $\dot k \ltimes k$ is the trivial one, the
action on $(\dot k / \dot k^N)$ is the one induced by the
$\mathbb Z / N \mathbb Z$-module structure of $\dot k / \dot k^N$
(see \cite{[Sc2]}), and the (right) actions on
Skew$_m(\mathbb Z / N \mathbb Z)$ and $(\mathbb Z / N \mathbb Z)^m$
are the natural ones (see also \cite{[Sc2]}).

Our description of the group structure of the biGalois group 
is different and seems to be simpler
than the one of Schauenburg in \cite{[Sc2]}, who used
both an action (of GL$_m(\mathbb Z / N \mathbb Z) \ltimes (\mathbb Z / N \mathbb Z)^m$) and a cocycle to describe the group structure. 
Let us try to understand the reasons leading to the use of a cocycle.
  
So let us consider a decomposable group datum $\mathbb G =(G,g,\chi)$
where $G$ is the direct product $G= \langle g \rangle \times K$.
We assume that $\chi^n = 1$ where $n = o(g)$.
An element of $\Gamma(\mathbb G)$ may be described as an element
$$
(f, \psi_1, \lambda , \overline{\sigma} , \psi_2)
\in {\rm Aut}(K) \times {\rm Hom}(K,\langle g \rangle) \times
\dot k \times H^2(K, \dot k) \times {\rm Hom}(K, \mu_n)$$
$${\rm satisfying} \quad
\chi(f(h)) \chi(\psi_1(h)) = \chi(h) \psi_2(h), \ \forall h \in K.$$
Hence we see that since $\chi^n = 1$, the element $\psi_2$ is redundant, and 
as a set we have 
$$\Gamma(\mathbb G) \cong 
{\rm Aut}(K) \times {\rm Hom}(K,\langle g \rangle) \times
\dot k \times H^2(K, \dot k).$$
The group law is not the one of a simple semi-direct product, and
this is where cocycles come into play for the
description of the group structure. However 
when $K = {\rm Ker} (\chi)$, then $\psi_1$ is redundant with respect
to the other data, and in this particular setting the
group law is the one of a semi-direct product, as described in Theorem 4.5.

\end{document}